%% file: pz1d.tex
%
%
\input amstex

\input eqnumams

\documentstyle{amsppt}
\overfullrule=0pt
\magnification =\magstep1
\baselineskip=18pt
\vcorrection{-.33truein}
\pageheight{9.0truein}
\input amssym.def

\def\fastvar{\nu}
\def\Epsilon{\Cal{V}}

\def\T3{{\tilde{\xi}}}
\def\BbbR{\Bbb{R}}

\def\BbbA{\Bbb{A}}

\def\BbbL{\Bbb{L}}
\def\BbbM{\Bbb{M}}

\def\CalA{\Cal{A}}

\def\CalE{\Cal{E}}

\def\CalL{\Cal{L}}
\def\CalO{\Cal{O}}
\def\CalQ{\Cal{Q}}
\def\CalR{\Cal{R}}

\def\bU{{\bar{u}}}

\def\R{\roman{Re\ }}

\def\sgn{\text{\rm sgn\ }}

\def\det{\text{\rm det\ }}

\def\Span{\text{\rm Span\ }}

\def\diag{\,\text{\rm diag}\,}

\def\Im{\text{\rm Im }}
\def\deltdot{{\cdot \atop {\raise8pt\hbox{$\delta$}}}}

\def\newsection{\centerline}


\pagewidth{6truein}
\pageheight{9truein}
\topmatter
\title
{An Evans function approach to spectral stability of small-amplitude
shock profiles
}
\endtitle
\leftheadtext{Stability of small amplitude shock profiles}
\rightheadtext{Ramon Plaza and Kevin Zumbrun}
\thanks
Research of the second author was supported
in part by the National Science Foundation under Grant No.  DMS-0070765.
\endthanks

\abstract
In recent work, the second author and various 
collaborators have shown using Evans function/refined 
semigroup techniques that, under very general circumstances,
the problems of determining 
one- or multi-dimensional nonlinear stability 
of a smooth shock profile may be reduced to 
that of determining spectral stability of the corresponding linearized 
operator about the wave.  
It is expected that this
condition should in general be analytically verifiable in the case of small
amplitude profiles, but this has so far been shown only on a case-by-case 
basis using clever (and difficult to generalize) energy estimates.  
Here, we describe how the same set of Evans function tools that were 
used to accomplish the original reduction can be used
to show also small-amplitude spectral stability by a direct and readily
generalizable procedure. This approach both recovers the results 
obtained by energy methods, and yields new results not previously obtainable.
In particular, we establish one-dimensional stability of small amplitude 
relaxation profiles, completing the Evans function program 
set out in Mascia\&Zumbrun [MZ.1].
Multidimensional stability of small amplitude viscous profiles
will be addressed in a companion paper [PZ], completing the
program of Zumbrun [Z.3].
\endabstract
\author
{Ramon Plaza and Kevin Zumbrun}
\endauthor
\date{April 5, 2002; Revised: June 2, 2002}
\enddate
\address
Courant Institute for Mathematical Sciences,
New York University,
215
Mercer Street, New York, NY 090010 (USA)
\endaddress
\email
plaza\@ cims.nyu.edu
\endemail 
\address
Department of Mathematics,
Indiana University,
Bloomington, IN  47405-4301
\endaddress
\email
kzumbrun\@indiana.edu
\endemail 
\endtopmatter
\document
\newsection {\bf Section 1. Introduction}
\sectionnumber=1 
\theoremnumber=0
\equationnumber=0
\smallskip
\TagsOnLeft

In this paper, we study one-dimensional
spectral stability in the small amplitude limit of smooth shock
profiles, i.e., traveling wave solutions 
$$
u=\bar u(x-st), \quad 
\lim_{z\to \pm \infty}\bar u(z)=u_\pm,
\eqnlbl{travelingwave}
$$
arising under various regularizations 
of a hyperbolic system of conservation laws
$$
u_t +  f(u)_{x}=0:
\eqnlbl{hyp}
$$
specifically, {\it viscous conservation laws} of form
$$
u_t + f(u)_{x}=(B(u) u_{x})_{x},
\eqnlbl{viscous}
$$
where $\R \sigma( B)\ge 0$,
and {\it relaxation systems} of form
$$
 \pmatrix u\\ v \endpmatrix_t
 +  \pmatrix \tilde f(u,v)\\ \tilde g(u,v) \endpmatrix_{x}
 = \pmatrix 0\\ q(u,v)
 \endpmatrix,
 \eqnlbl{relaxation}
$$
where 
$$
\text{\rm Re }\sigma\big(q_v(u,v^*(u))\big)<0
\eqnlbl{qv}
$$
along a smooth {equilibrium manifold} defined by 
$$
q(u,v^*(u))\equiv 0,
\eqnlbl{eq}
$$
and $f(u)=\tilde f(u,v^*(u))$.
Here, $x$, $t\in \BbbR^1$, $u$, $f$, $\tilde f \in \Bbb{R}^n$, $v$, $\tilde g$,
$q \in \Bbb{R}^r$,
and $B\in \BbbR^{n\times n}$, where typically $n$, $r>1$;
we make the standard
assumptions of strict hyperbolicity of $Df$,
and genuine nonlinearity of the principal characteristic eigenvalue (the one
associated with the approximate shock direction), along with further standard
hypotheses on $B$ and $\tilde f$, $\tilde g$, $q$, 
to be described later on: in particular,
that systems \eqnref{viscous} and \eqnref{relaxation} be {dissipative}
in the sense of Kawashima and Zeng [Kaw, Ze].
These are generically satisfied for 
most models of physical interest;
for discussions of applications, we refer the reader to the
general surveys [Z.3,Z.6,MZ.1,N,Yo.4].

A recent development in the stability analysis of smooth shock
profiles has been the successful adaptation of Evans function/dynamical
system techniques
to this nonstandard setting, with an associated
explosion of new results; see, e.g.,
[GZ,H.1--3,ZH,BSZ,BMSZ,ZS,HZ,Z.2--4,HoZ.1--2,OZ.1--2,MZ.1--2].
For the origins of these methods
in the setting of reaction diffusion and singular perturbation problems,
see [E.1--4,J,AGJ,PW]; see also the important analyses of [Sat,K.1--2,JGK] in
the scalar case.
In particular, it has now been established under extremely
general circumstances that both linearized and nonlinear stability
are implied by {spectral stability} of the linearized operator
about the profile, a fact that was left in doubt by earlier
``direct'' analyses [MN,Go.1--3,KMN,SX.1--2,S,GoM,LZ.1--2,FreL,L.1--3] 
carried out on a case-by-case basis.
This reduces the question of stability to an ODE issue, and raises
the hope that at least small amplitude stability can be treated in
a uniform way across dimension and type of regularization.

On the other hand, spectral stability has up to now been
verified only by energy estimates resembling those of the earliest
direct analyses of [MN,Go.1--2,KMN]; see [Z.1,HuZ,Liu,Hu.1,Z.4] for
examples in various contexts.
And, these have proven difficult to generalize: in particular,
one-dimensional stability of general small amplitude relaxation profiles, and
multi-dimensional stability of small amplitude viscous shock profiles remain
outstanding open problems in the theory.
In this, and a companion paper [PZ] treating the multidimensional case,
we remove this troublesome gap in the theory, showing that the
same Evans function tools that were used to effect the reduction to
spectral stability may be used as the basis of a general procedure
to verify spectral stability in the zero-shock-amplitude limit,
applicable in particular in the two open cases mentioned just above.
This both unifies and completes existing stability theory, at least as 
regards the small amplitude case; 
stability of large amplitude shock profiles remains as the outstanding
open question in this area.

\medskip
The basic idea behind our approach is a natural one suggested
by Gardner and Jones [GJ.3] in the context of the one-dimensional
strictly parabolic case.
In the weighted energy method of Goodman [Go.1--2], the
approach on which most subsequent extensions have been based,
the strict transversality of characteristic fields other than
the principal one is used to ``project out'' behavior in transverse
modes, reducing the problem to an approximate scalar conservation 
law in the principal characteristic field.
Gardner and Jones pointed out that transverse modes correspond
in the associated eigenvalue ODE roughly to ``fast'' and 
``super-slow'' modes, which on a large portion of the relevant
spectral domain may be projected out using dynamical system
techniques to leave an approximately scalar ``slow'' manifold flow
as in the argument of Goodman, and proposed this observation
as the basis of a dynamical systems argument alternative
to the one-dimensional stability argument of Goodman.
However, Gardner and Jones did not examine the crucial small-frequency
regime where super-slow and slow modes intermingle,
and therefore  eliminated only the possibility of eigenvalues with real part 
greater than 
a small, but fixed constant $\theta>0$;
by contrast, a standard G\"arding-type energy
estimate gives the much stronger result $\R \lambda \le C\epsilon^2$,
$\Im \lambda \le C\epsilon$,
where $\epsilon:=|u_+-u_-|$ denotes shock amplitude.
Thus, at a technical level, substantial issues still remain;
indeed, virtually all of our effort will be directed to the understanding
of the ``mixing'' regime complementary to that considered in [GJ.3].

\medskip
In the remainder of the introduction, we give an overview of the analysis
and main results,
deferring technical aspects to the main body of the paper.

\medskip
{\bf Spectral stability.} Let us first make precise the notion of
spectral stability that we consider here.  
Without loss of generality taking $s=0$
in \eqnref{travelingwave} (by shifting to an appropriate 
traveling frame), we have that $\bar u(x)$ is a stationary
solution of the nonlinear evolution system 
$$
u_t=\Cal{F}(u)
\eqnlbl{abev}
$$
described by its associated regularized conservation law,
of form \eqnref{viscous} or \eqnref{relaxation}.
Linearizing this system about the stationary solution $\bar u$,
we obtain a linearized evolution system
$$
v_t=Lv,
\eqnlbl{ablinear}
$$
where $L:=d\Cal{F}(\bar u)$ is the formal derivative of operator $\Cal{F}$
about $\bar u$.

{\bf Definition \thmlbl{specstab}.} Following [ZH,Z.3],
we define strong spectral stability (condition (D1') of [Z.3]) as
$$
\sigma_p (L) \subset \{\lambda: \R \lambda <0\} \cup \{0\}.
\eqnlbl{1dspecstab}
$$
where $\sigma_p(M)$ denotes point spectrum of a linear operator $M$:
equivalently, there exist no $L^2$ solutions of the eigenvalue
equation $(L-\lambda)w=0$ for $\R\lambda \ge 0$ and $\lambda \ne 0$.
\medskip

{\bf Remark \thmlbl{nonlinstab}.}
It has been shown in a variety of contexts, 
in particular, all contexts considered here,
that strong spectral stability together with
the ``hyperbolic'' conditions of inviscid stability of the corresponding ideal
shock $(u_-,u_+)$ of \eqnref{hyp} (dynamical, or ``outer'' stability)
and transversality of the traveling wave connection (structural,
or ``inner'' stability) implies linearized and nonlinear
stability; see [ZH,Z.2,Z.3,MZ.1--2].
In the small amplitude case, under reasonable assumptions,
inviscid stability holds always;
see [M.1--3] in the one-dimensional case, [M\'e.1--5,FM,FZ,Z.5--6]
in the multidimensional case.  Likewise, under reasonable assumptions,
connections are always transverse; see discussion of existence theory
given below.
Thus, {\it strong spectral stability is sufficient to imply linearized
and nonlinear stability}. 
\medskip

\medskip
{\bf The canonical model.}
In the one-dimensional genuinely nonlinear context, 
it is well known that {\it Burgers equation},
$$
u_t + (u^2/2)_x= u_{xx},
\eqnlbl{burgers}
$$
approximately describes small-amplitude viscous behavior in the principal 
characteristic mode.
In particular, the family of exact solutions
$$
\bar u^\epsilon(x):= -\epsilon \tanh(\epsilon x/2)
\eqnlbl{burgersprofile}
$$
give an asymptotic description of the structure of weak
viscous shock profiles in the principal direction, in the
limit as amplitude $\epsilon$ goes to zero.
Note that \eqnref{burgers} is invariant under the parabolic
scaling $x\to \epsilon x$, $t\to \epsilon^2 t $, $u\to u/\epsilon $,
whence we may immediately conclude that any unstable spectrum of
the linearized operator about the profile $\bar u^\epsilon$ must
lie within a ball of radius $C\epsilon^2$.
(By standard considerations,
the unstable spectrum of a second order parabolic operator 
is at least bounded).

Burgers profiles are strongly spectrally stable, as
may be readily established by a 
number of different techniques:
for example, by $L^1$ contraction [Kr], Sturm-Liouville 
considerations [Sat,He,Z.5],
or direct energy estimates as in [MN,Go.1--2,Z.3 Appendix A.6)].
We record this fact as:

\proclaim{Proposition \thmlbl{modelstab}}
Shock profiles \eqnref{burgersprofile} (of any amplitude $\epsilon$,
by scale-invariance) 
are strongly spectrally stable as solutions of \eqnref{burgers}.
\endproclaim

\medskip
{\bf The reduced profile equations.}
Existence/structure of small-amplitude profiles
may be determined by center manifold reduction,
as pointed out by Majda and Pego [MP]
in their pioneering treatment of the general strictly parabolic case;
extensions to general partially parabolic, or ``real'' viscosity,
and relaxation systems have been carried out in [P], [Yo.Z,FZe,MZ.1],
respectively.
The result in all cases, is that the center manifold associated with
the weak profile problem (appropriately chosen) is foliated by one-dimensional
fibers, on which the reduced flow, in the genuinely nonlinear case, is
approximately that of the profile equation for 
Burgers equation \eqnref{burgers}.

More precisely, the center manifold lies tangent to (i.e., 
within angle $\Cal{O}(\epsilon)$ of) the principal eigendirection $r_p$
of $Df^1$, and the reduced flow, parametrized by 
$\eta:=l_p(u_-)\cdot (\bar u- (u_- + u_+)/2)$,
$l_p$ the corresponding left eigendirection, obeys the approximate
Burgers profile equation
$$
\beta\eta'= \frac{\Lambda}{2}(\eta^2- \epsilon^2) + \Cal{O}(|\eta,\epsilon |^3),
$$
or, rescaling by 
$$
\eta\to \eta/\epsilon, \quad 
x\to \Lambda \epsilon x/\beta:
\eqnlbl{rescaling}
$$
$$
\eta'= \frac{1}{2}(\eta^2- 1) + \epsilon O(\eta,\epsilon),
\eqnlbl{rescaledprofile}
$$
$O(\cdot,\cdot)\in C^1$,
as compared to the exact profile equation
$$
\bar \eta'= \frac{1}{2}(\bar \eta^2- 1)
\eqnlbl{exactprofile}
$$
for \eqnref{burgers}--\eqnref{burgersprofile} with $\epsilon=1$;
here and below, ``$'$'' denotes $\partial/\partial x$.
The effective diffusion coefficient $\beta$ is given by
$\beta:=l_p Br_p(u_-)$
in the viscous case, and
$\beta:=l_p B_*r_p(u_-)$ in the relaxation case, where 
 $$
  \eqalign{
   B_*(u)
    &:= -\tilde f_v q_v^{-1}(g_u - v^*_u f_u) \cr
    &= -\tilde f_v q_v^{-1}(\tilde g_u - \tilde g_vq_v^{-1}q_u
     + q_v^{-1}q_u(\tilde f_u - \tilde f_v q_v^{-1} q_u), \cr
           }
  \eqnlbl{Bstar}
 $$ 
$g:=\tilde g(u,v^*(u))$, 
denotes the effective ``Chapman--Enskog''
viscosity associated with the relaxation system,
and $\Lambda:=l_p^t D^2f^1(r_p,r_p)(u_-)\ne 0$ denotes the coefficient
of genuine nonlinearity for the hyperbolic system \eqnref{hyp}; 
for details, see [MP,Z.1,P,Yo.Z,FZe,MZ.1].

From this reduction to normal form we obtain in standard fashion:

\proclaim{Proposition \thmlbl{limprofile}}
Under the rescaling \eqnref{rescaling}, there hold
$$
|\eta_\pm - \bar \eta_\pm |\le C\epsilon
\eqnlbl{iprof}
$$
and
$$
|(\eta- \eta_\pm)-(\bar \eta -\bar\eta_\pm)|\le C\epsilon e^{-\theta |x|}
\eqnlbl{iiprof}
$$
for $x\gtrless 0$ for any fixed $0<\theta<1$, and some $C=C(\theta)>0$,
where $\eta_\pm$ and $\bar \eta_\pm:=\pm 1$ denote the
rest points of \eqnref{rescaledprofile} and \eqnref{exactprofile},
respectively.
\endproclaim

{\bf Proof.}  Assertion \eqnref{iprof} follows by the Implicit
Function Theorem from normal hyperbolicity of $\bar\eta_\pm$.
Assertion \eqnref{iiprof} follows from standard 
stable/ustable manifolds estimates.
\qed
\medskip
\proclaim{Corollary\thmlbl{limchar}}
Under the rescaling \eqnref{rescaling}, the principle
eigenvalue $a_p:= a_p(\bar u)$ converges
to $\bar \eta$ at the same rates \eqnref{iprof}--\eqnref{iiprof} 
as does $\eta$.
\endproclaim

{\bf Proof.} This follows using tangency of the profile to direction $r_p$
together with the fact [Sm] that $\nabla a_p \cdot r_p= \Lambda$.
\qed

\medskip
{\bf The reduced eigenvalue equations.}
Following [GJ.3], we carry out here a reduction of the
generalized eigenvalue equations similar in spirit to
that carried out by Majda\&Pego and others for the profile equations.
However, whereas those anlayses involved center manifold reduction of
an autonomous nonlinear system, 
we shall work directly with the linear, nonautonomous system
given by the eigenvalue equation $(L-\lambda)w=0$
associated with \eqnref{ablinear}, written as an appropriate
first-order system
$$
W'=\BbbA(x,\lambda) W,
\eqnlbl{firstorder}
$$
using the {\it tracking lemma} of [GZ,ZH] (specifically,
a refined version introduced in [MZ.1])
to effect the reduction to a slow manifold; see Section 2.
Here, \eqnref{firstorder} is dimension $2n$ in the strictly
parabolic viscous case, $n+ \text{\rm Rank} B$ in the
general viscous case, and $n+r$ in the relaxation case;
as above, ``$'$'' denotes $\partial/\partial x$.
For the efficient coordinatization of the first-order system 
\eqnref{firstorder} in the relaxation and general viscous case
(important for practical computation), we follow the scheme of [Z.3];
see sections 3--4 for details.

Of course the tracking lemma is itself an analytic formulation of
classical invariant manifold techniques; 
see the original formulation in [GZ]
in terms of projectivized flow/invariant sets of autonomous
systems.
The calculation of the reduced flow on the slow manifold can equally well
be carried out by classical center manifold reduction;
see the concurrent treatment of Plaza [Pl] in the one-dimensional case;
see also the independent treatment of [FSz] (described further in
the note following the introduction). 
Indeed, this may be preferable in the more complicated situations 
arising with various types of degeneracies in the underlying equations.
However, there is some advantage in taking account of the linear
nature of the underlying problem, in that the calculation of reduced
equations simplifies considerably.

\medskip

{\it First reduction.}
A first application of the tracking lemma projects out ``fast''
transverse modes: $n-1$ in the strictly parabolic viscous case,
$\text {\rm Rank}B-1$ in the general viscous case, and
$r-1$ in the relaxation case, in each case reducing to an 
$(n+1)$-dimensional ``slow manifold,''
of which $n-1$ correspond to ``super slow'' transverse modes $\rho_j$, 
and the remaining $2$ to ``medium slow'' Burgers modes which are 
roughly linearized versions of $\eta$, $\eta'=:z$ 
in the limiting equations above.
Rescaling by
$x\to \Lambda\epsilon x/\beta$, 
$\lambda\to \beta\lambda/\Lambda^2\epsilon^2$,
and (balancing by) $\eta'\to \beta \eta'/\Lambda\epsilon$, 
we obtain a normal form for the reduced flow that has
bounded ($\Cal{O}(1)$) coefficients for bounded $\lambda$
and, modulo error terms  $\Phi_j(x,\lambda) \tilde W$ 
($\tilde W:=(\rho,\eta,z)$),
becomes a block triangular system in $\rho=(\rho_+,\rho_-)$ 
and $(\eta,z)$, of form
$$
\rho_\pm'= M_\pm(x,\lambda)\rho_\pm,
\eqnlbl{rhoeqn}
$$
and
$$
\pmatrix \eta\\ z\endpmatrix'=
\bar M_0(x,\lambda)\pmatrix \eta\\ z\endpmatrix + N(x,\lambda)\rho,
\eqnlbl{etaeqn}
$$
where: 
$$
\Phi_{\rho_\pm}(\pm \infty, \lambda)\equiv 0, \qquad 
\Phi_{(\eta,z)}(\pm \infty,\lambda)= 
\epsilon
\Cal{E}_\pm(\lambda,\epsilon)
\eqnlbl{inftybds}
$$
$$
|\Phi(x,\lambda)-\Phi(\pm\infty, \lambda)|\le 
C(1+|\lambda|)\epsilon e^{-\theta |x|}
\quad
\text{\rm for $x\gtrless 0$},
\eqnlbl{Phibds}
$$
and
$$
|N|\le Ce^{-\theta |x|};
\eqnlbl{Nbds}
$$
$$| M_\pm|=\Cal{O}(\epsilon|\lambda|), \qquad
M_\pm \gtrless \epsilon\theta(|\R \lambda|+ \epsilon^2|\Im \lambda|^2), 
\eqnlbl{superslow}
$$
$\theta>0$; and $\bar M_0$ in the viscous case corresponds exactly
to the coefficient 
$$
\bar \BbbA_0:=
\pmatrix 
0& 1 \\
\lambda & \bar \eta(x)
\endpmatrix
\eqnlbl{burgerseval}
$$
of the eigenvalue equation 
for the standard Burgers equation 
\eqnref{burgers}--\eqnref{burgersprofile},
$\epsilon=1$,
written as a first-order system.
Here, $\Cal{E}_\pm(\cdot)$ is a smooth function satisfying
the growth bounds
$$
\Cal{E}_\pm=
\pmatrix \Cal{O}(1 +|\lambda| \epsilon) & \Cal{O}(1+|\lambda|\epsilon)\\
\Cal{O}(1+ |\lambda|)  & \Cal{O}(1+|\lambda|)\\
\endpmatrix.
\eqnlbl{Einftybds}
$$
{\bf Remark \thmlbl{needed}.} The asymmetric bounds on $\CalE_\pm$
are a result of rescaling.  The stronger bound
in the upper righthand corner is actually
important for our later analysis, specifically, in the
large-frequency regime where we approximately undo 
the rescaling in the course of our argument.

\medskip
{\it Further reductions/normal forms.}
In the rescaled coordinates, there are three distinct frequency regimes:
\smallskip
I. $|\lambda|\le C$;
II. $C\le |\lambda|\le C\epsilon^{-1}$;
and III. $C\epsilon^{-1} \le |\lambda|\le C^{-1}\epsilon^{-2}$, 
\smallskip
\noindent
$C>0$ sufficiently large,
corresponding in the original (unrescaled) coordinates to
$|\lambda|\le C\epsilon^2$,
$C\epsilon^2 \le |\lambda|\le C\epsilon$,
and $C\epsilon \le |\lambda|\le C^{-1}$, respectively.
(By standard considerations, the spectrum of $L$ is
restricted to small frequencies $|\lambda|\le C^{-1}$;
see Proposition 2.9, below.)
Each regime requires a slightly different treatment, and is associated
with a different normal form.

\smallskip
Region III, the region considered by Gardner and Jones [GJ], we will
denote as the ``parabolic'' regime, as this is the regime on which
behavior is dominated by dissipativity of the underlying
system (which yields ``parabolic'' behavior for small frequencies,
as evidenced by the structure of the reduced system
\eqnref{rhoeqn}--\eqnref{etaeqn}),
independent of any other structure of the traveling wave profile whatsoever.
This regime may be treated by the standard ``high-frequency'' techniques
of the usual Evans function theory (see [AGJ,GZ,ZH,MZ.1]), which
are roughly equivalent to sectorial, or G\"arding-type, energy estimates,
to preclude altogether the possibility of eigenvalues.

Specifically, $(\eta,z)$ may be resolved into decaying/growing
modes $\eta_\pm$ that have a uniform spectral gap both from each other
and from $\rho_+$ and $\rho_-$, the latter of which also have
a uniform gap from one another, to an order sufficient that a
second application of the 
tracking lemma, to the once-reduced equations \eqnref{rhoeqn}--\eqnref{etaeqn}
yields further reductions of the flow onto four transverse 
invariant manifolds tangent to each of the $\rho_\pm$, $\eta_\pm$ directions.  
In particular, we find that the manifold of solutions decaying at $+\infty$ 
is transverse to that of solutions decaying at $-\infty$, yielding the result.
In other words, the normal form for this regime is trivial, with growing
and decaying modes decoupled.

Region II we denote as the ``reduced parabolic'' regime, since this
is the regime for which similar tracking/energy estimates prohibit
eigenvalues for the Burgers equation.
Here also, we find that eigenvalues cannot exist, independent of the
structure of the profile.  However, the analysis is much more delicate;
indeed, we regard this as the trickiest case in our argument.
Here, we can again separate off decay/growth modes $\eta_\pm$ (parabolic
behavior of the Burgers part) as in the previous case, but $\rho_\pm$
are not sufficiently separated to apply the tracking lemma to these
coordinates.  In this case, therefore, our ultimate reduction is to
the $(n-1)$-dimensional superslow equations \eqnref{rhoeqn}, but with
error terms now involving $\rho$ alone and not $(\eta,z)$.
Note that the coefficients $M_\pm$ are bounded on regime II, so that
this is indeed a valid normal form.

Region I we denote as the ``gap regime,'' in reference to the fact
that on this regime we rely solely on the {\it gap lemma} of [GZ,KS,ZH,Z.3]
(specifically a refined version given in [M\'eZ]) in our analysis
of the slow flow, rather than
using the tracking lemma to carry out a further reduction as in
Regions II and III.
The ``gap'' in ``gap lemma'' refers in fact to absence of spectral
gap, as is the case in this regime; the lemma asserts 
that {\it uniform exponential convergence} of 
coefficients as $x\to \pm \infty$
and {\it continuous extension} of the stable/unstable eigenspaces of
the limiting coefficient matrices $\BbbA_\pm(\lambda)$ at $x\to \pm \infty$
can substitute in this situation for uniform spectral gap, to yield
estimates valid on $x\gtrless 0$.
In this regime, which is the crucial one determining behavior,
the normal form is the entire $(n+1)$-dimensional system 
\eqnref{rhoeqn}--\eqnref{etaeqn}.

\medskip
{\bf The gap lemma: convergence to limiting flows.}
We complete our reduction in each case by an application of the gap lemma
(specifically, a refinement given in [M\'eZ]),
showing that the flow of the normal forms we have obtained
converges to that of the ``limiting'' normal forms obtained
by dropping error terms $\Theta \tilde W$, etc.: more precisely,
the flows associated with the stable manifold at $+\infty$
and the unstable manifold at $-\infty$, or in some cases
their analytic continuation into the essential spectrum.
The basic approach is to: (i) 
conjugate approximate and limiting
equations on the half-lines $x\gtrless 0$ to their asymptotic 
systems at $\pm \infty$ by a linear change of variables uniformly
close to the identity
(an application of the gap lemma,
made possible by the uniform exponential decay estimate \eqnref{inftybds}), 
then (ii) check by direct
linear algebraic computation that the stable and unstable
subspaces of the respective asymptotic coefficient matrices are uniformly close.
In the present applications, (ii) is straightforward, since
the asymptotic matrices of approximate and limiting flows agree
for the $\rho$ equation, while for the 
(decoupled) $(\eta,z)$ equation, the
stable and unstable subspaces for the limiting equations are, under
our strategically chosen rescaling, uniformly spectrally separated
on the region I where it comes into play.

Convergence of the stable flow at $+\infty$ and unstable flow at
$-\infty$ then implies uniform convergence as $\epsilon \to 0$
of suitably chosen {\it Evans functions} associated with the approximate
and limiting normal forms.
The Evans function is a Wronskian measuring solid angle between the
stable and unstable subspaces of a given eigenvalue ODE written
as a first-order system, constructed so as to be analytic in the frequency
$\lambda$ on a ``region of consistent splitting''
(defined (A1), Section 2).
In each of the contexts considered here, the region of consistent splitting
includes the entire set of frequencies
$\{ \R \lambda \ge 0 \}\setminus \{0\}$
of interest; see, e.g., [Z.3].  Moreover, both the Evans function and 
its component subspaces {\it extend continuously along rays} in frequency
space, in both the one- and multidimensional case [ZS,Z.3] to
the closure at the origin, making possible uniform estimates
of transversality; indeed, in the one-dimensional case considered here,
the Evans function may in fact be extended analytically through
the origin, without restricting to rays [GZ].
(Note: in the multidimensional case, the Evans function
may be holomorphically extended along rays [ZS,Z.3]; however, in 
certain ``glancing'' directions, there arise branch singularities 
at the origin.)
On the region of consistent splitting,
the zeroes of the Evans function correspond both in location and multiplicity 
with the eigenvalues of the associated linear operator;
for details, see [GJ.1--2,ZH].

The multiplicity of the generalized eigenvalue at the origin
$\lambda=0$ is directly calculable and related
to hyperbolic stability, as shown in the one- and 
multi-dimensional case in [GZ,BSZ] and [ZS], respectively.
In the present, small-amplitude case, it is readily seen to
be the same as the multiplicity for the limiting normal form
equations, namely one. 
In the one-dimensional, analytically extendable case, 
we may thus conclude immediately from uniform convergence to the limiting
Evans function (by degree, i.e., winding number, considerations) 
that {\it strong spectral stability
of the original system is equivalent to strong spectral
stability for the limiting normal forms}.
Alternatively, we may remove the zero of the Evans function
at the origin by working instead with the ``integrated'' eigenvalue
equations following [MN,Go.1] or with ``flux variables'' following [Go.2],
to recover the same conclusion using only continuity at the origin of $D$.

In this paper, we follow the latter strategy exclusively.
In the multidimensional analysis of [PZ], we use a 
modification of the same strategy,
specifically, an interesting variation intermediate to 
flux and integrated variable methods, neither of which 
themselves directly generalize to multidimensions in a useful way,
to achieve the same result;  this ``balanced flux form'' 
is quite similar to that introduced in the
treatment of relaxation systems in Section 4 of the present paper.
The former strategy is perhaps viable as well, using 
computation of winding number on Riemann surfaces as suggested in [GZ];
however, this would involve extra bookkeeping (specifically,
tracking of branch points in approximate vs. limiting systems) 
that we prefer to avoid.

\medskip
{\bf Conclusions/summary of the main results.}
From this point, we may immediately deduce our main results.
For, in each of the (nontrivial) limiting normal forms
arising in regions I--II, the equations for
$\rho_\pm$ completely decouple, and clearly support
no ``decaying'' solutions other than the trivial, zero solution
(in general, ``decaying'' is defined as lying in $\rho_+$
direction at $+\infty$ and $\rho_-$ direction at $-\infty$;
away from the essential spectrum, this is equivalent to actual
decay at $\pm \infty$).
Thus, we may conclude that $\rho\equiv 0$ for any generalized
eigensolution, leaving a decoupled $(\eta,z)$
equation in region I, {\it agreeing with
the generalized eigenvalue equations for
canonical model \eqnref{burgers}},
and the trivial flow in region II.
Recall that the reduced flow was already trivial in region III.
Likewise, $\fastvar\equiv 0$ in each of these regions,
since fast growing and fast decaying modes completely decouple.

Noting that the canonical model, being among the class considered,
is subject to the same reduction in regions II--III,
we may conclude, finally, our main result: 
that {\it strong spectral stability is equivalent to strong spectral
stability of the canonical model}, at least modulo the degenerate
case that the canonical model yields nontrivial nonstable eigenvalues 
lying precisely on the imaginary axis; see Proposition 3.2 and
Remark 3.3 below.
As a corollary, recalling the result of Proposition \thmref{modelstab}, 
we obtain the stated results of strong spectral stability of
small-amplitude profiles, across the general class of models considered
in \eqnref{viscous} and \eqnref{relaxation}.
More precisely, in this paper, we carry out the details of the 
one-dimensional case for general strictly parabolic viscosities,
and for general relaxation systems.  The case of general real
viscosity may be carried out similarly as the relaxation case,
following the dual framework set out in Appendix A.2 of [Z.3];
we omit these calculations, 
as satisfactory one-dimensional results already exist [HuZ].

\bigskip
{\bf Plan of the paper.}  In Section 2, we recall and
slightly extend the basic Evans function tools we will need.
In Section 3, we carry out the one-dimensional, strictly parabolic
viscous case.  For clarity, we first carry
out in its entirety the simpler identity viscosity case,
afterward indicating the necessary adjustments in the general case.
In Section 4, we carry out the general relaxation case, for
clarity first treating the $2\times 2$ case by reduction
to the scalar viscous case; we note that this gives an independent 
proof of stability in this case, different from that of
Liu [L.2].
Finally, in Section 5, we briefly indicate the extension of our analysis
to the multidimensional viscous case; details will be given in [PZ].

\bigskip
{\bf Note}:  
Shortly before the completion of this analysis
(precisely: after the completion of the one-dimensional
viscous case and before the completion of the relaxation case),
we have learned of closely related work of
H. Freist\"uhler and P. Szmolyan [FreS] 
in which they
establish a reduction method similar (indeed, perhaps equivalent)
to ours, but proceeding entirely from 
the point of view of geometric singular perturbation theory.
In the cited work, the first in a planned series of
three, they carry out a complete analysis of the 
one-dimensional identity viscosity case,  announcing the intention
to treat one-dimensional relaxation/real viscosity and multidimensional
stability in papers two and three, respectively.

\bigskip
\newsection {\bf Section 2. Evans function framework.}
\sectionnumber=2 
\theoremnumber=0
\equationnumber=0
\smallskip
\TagsOnLeft

We begin our analysis by recalling, and in some cases slightly
extending, the Evans function tools we will need.
For later reference, we carry out the analysis in generality
sufficient for the multidimensional case as well.

\medskip
{\bf 2.1. The gap lemma and convergence of approximate flows.}
Consider a family of first-order systems
$$
{W}'=\BbbA^\epsilon(x,\tilde \xi, \lambda) W
\eqnlbl{gfirstorder}
$$
indexed by parameter $\epsilon$,
where $(\tilde \xi,\lambda) $ vary within the set
$$
\Omega:=\{(\tilde \xi, \lambda): \,
\tilde \xi \in \BbbR^{d-1}, \, \R \lambda \ge 0\}
\eqnlbl{omega}
$$
of unstable--neutrally stable frequencies, $\epsilon$ varies
within some given open set $\Epsilon$, and $x $ varies within $\BbbR^1$.
Equations \eqnref{gfirstorder} are to be thought of as generalized eigenvalue
equations, with parameter $\lambda$ representing a Laplace transform
frequency in time, and $\tilde \xi$ representing 
a Fourier transform frequency  in directions
of spatial symmetry; in the present, one-dimensional case, 
$\tilde \xi\equiv 0$.
In the applications of this paper and of [PZ], $\epsilon$ will be just the
shock strength $|u_+-u_-|$. 

We make the following basic assumptions.

\medskip

(A0) \quad 
Coefficient $\BbbA^\epsilon(\cdot,\tilde \xi,\lambda)$, considered
as a function from $(\tilde \xi, \lambda,\epsilon)$ into $L^\infty(x)$
is analytic in $(\tilde \xi,\lambda)$ and $C^k$ in $\epsilon$ for
some $k\ge 0$ on $\Omega \times \Epsilon$.
Moreover, $\BbbA^\epsilon(\cdot, \tilde \xi,\lambda)$ approach
exponentially to limits $\BbbA_\pm$ as $x\to \pm \infty$, 
with uniform exponential decay estimates
$$
|\BbbA^\epsilon- \BbbA^\epsilon_\pm| \le C_1e^{-|x|/C_2}, \, C_j>0, \quad
\text{\rm for } x\gtrless 0
\eqnlbl{expdecay}
$$
on compact subsets of $\Omega \times \Epsilon$.
\smallskip
(A1) \quad 
 On $(\Omega \setminus \{(0,0)\})\times \Epsilon$, the 
limiting matrices $\BbbA^\epsilon_+$ and $\BbbA^\epsilon_-$ 
are both hyperbolic (have no center subspace), and the 
dimensions of their stable (resp. unstable)
subspaces $S^\epsilon_+$ and $S^\epsilon_-$ 
(resp. unstable subspaces $U^\epsilon_+$ and $U^\epsilon_-$) agree.
\smallskip
(A2) \quad 
At the origin $(\tilde \xi,\lambda)=(0,0)$,
subpaces $S^\epsilon_\pm$, $U^\epsilon_\pm$ have continuous (to some order,
which may be even analytic, or holomorphic with branch point at the origin) 
limits along rays (i.e.,  for
$(\tilde \xi,\lambda)=(r \tilde \xi_0,r \lambda_0)$ as $r\in \BbbR \to 0^+$,
$(\tilde \xi_0, \lambda_0)\in \Omega\setminus \{(0,0)\}$),
which, moreover, are $C^k$ in the parameter $\epsilon$, where $k$
is as in (A0) above.
\medskip

These are satisfied in all the contexts considered in 
both this paper and [PZ]; see [Z.3].
Condition (A0) is induced, at least for $\epsilon$ bounded from 
zero,  by the origins of \eqnref{gfirstorder}
through the linearization about smooth shock profiles of
\eqnref{viscous} or \eqnref{relaxation}.
Uniform approach \eqnref{expdecay} as $\epsilon\to 0$ will
be recovered through reduction/rescaling.
Condition (A1) may be recognized as the standard hypothesis
of ``consistent splitting,'' as introduced in [AGJ],
and follows by the assumption of dissipativity of systems
\eqnref{viscous} and \eqnref{relaxation}.
Condition (A2) is an extension suitable for multidimensions
of the ``geometric separation'' hypothesis of [GZ]
and is a consequence of the hyperbolic structure of \eqnref{hyp}.
(Note: separation was not needed for the arguments of [GZ], 
and in fact does not hold in multidimensions; see discussion of
Appendix A.4, [Z.3].) 

\medskip

{\bf The gap lemma.} 
Under these circumstances, the ``gap lemma'' established
in various degrees of generality in [GZ,KS,ZH,Z.3]
asserts that behavior of the variable-coefficient equation
\eqnref{gfirstorder} may be related to that of its constant-coefficient
limiting systems on $x\gtrless 0$, while maintaining the assumed
regularity in all parameters.
This has recently been greatly improved in [M\'eZ], in the form of the
following {\it conjugation lemma}, a version that is particularly
convenient in applications.

\proclaim{Lemma \thmlbl{conjugation}[M\'eZ]}
Under assumptions (A0)--(A2), there exists locally to any 
given $(\tilde \xi, \lambda, \epsilon)\in \Omega \times \Epsilon$
a pair of linear transformations
$P^\epsilon_+(x,\tilde \xi,\lambda)$ and $P^\epsilon_-(x,\tilde \xi,\lambda)$
on $x\ge 0$ and $x\le 0$, respectively,
and possessing the same
regularity as $\BbbA^\epsilon$ in all arguments, such that:
\medskip
(i)
$|P^\epsilon_\pm|$ and their inverses are uniformly bounded, with
$$
|P^\epsilon_\pm - I|\le C C_1 C_2 e^{-\theta |x|/C_2}
\eqnlbl{Pdecay} 
$$
where $0<\theta<1$ is an arbitrary fixed parameter, and $C>0$ 
and the size of the neighborhood of definition depend
only on $\theta$,  the modulus of the entries of 
$\BbbA^\epsilon$ at $(\tilde \xi, \lambda, \epsilon)\in \Omega \times \Epsilon$,
and the modulus of continuity of $\BbbA^\epsilon$ on some neighborhood
of $(\tilde \xi, \lambda, \epsilon)\in \Omega \times \Epsilon$.
\smallskip
(ii)  The change of coordinates $W:=P^\epsilon_\pm Z$ reduces \eqnref{gfirstorder}
on $x\ge 0$ and $x\le 0$, respectively, 
to the asymptotic constant-coefficient equations
$$
Z'=\BbbA^\epsilon_\pm(\tilde \xi, \lambda) Z.
\eqnlbl{glimit}
$$
\endproclaim
\medskip

{\bf Proof.} 
As described in [M\'eZ],
this is a straighforward corollary of the gap lemma 
as stated in [Z.3], applied to the ``lifted'' matrix-valued equations 
for the conjugating matrices $P^\epsilon_\pm$.
\qed

\medskip
{\bf Definition \thmlbl{evans}.}
Let $v_{1+}^\epsilon, \dots, v_{k+}^\epsilon$ and 
$v_{(k+1)-}^\epsilon, \dots, v_{N-}^\epsilon$
be bases for $S^\epsilon_+$ and $U^\epsilon_-$, as defined in (A1), chosen
with the same regularity assumed on $S^\epsilon_+$, $U^\epsilon_-$
(note: that such a choice is possible is a consequence of
standard matrix perturbation theory [Kat]),
and $P^\epsilon_\pm$ be transformations defined as above
on some neighborhood of $(\tilde \xi,\lambda,\epsilon\in \Omega\times \Epsilon$.
Then, the {\it local Evans function} for \eqnref{gfirstorder}
associated with these choices is defined as
$$
D^\epsilon(\tilde \xi, \lambda):=
\det
\pmatrix 
P^\epsilon_+ v_{1+}^\epsilon& \dots& P^\epsilon_+ v_{k+}^\epsilon&
P^\epsilon_- v_{(k+1)-}^\epsilon& \dots&P^\epsilon_- v_{N-}^\epsilon
\endpmatrix_{|x=0}.
\eqnlbl{evans}
$$

\medskip
{\bf Remark \thmlbl{nonunique}.}
Evidently, the choice of local Evans function is highly nonunique.
Though we shall not need it here,
it can be shown using the original formulation of the gap lemma
in [GZ,KS] that a globally analytic choice is possible in both
the one- and multidimensional case; see [GZ] and [ZS,Z.3], respectively.

\medskip
Combining \eqnref{glimit} and \eqnref{evans}, we immediately obtain
the following result sufficient for our needs.
\footnote{
Recall that Gardner and Jones [GJ.1--2] have established the
stronger result of correspondence up to multiplicity between
eigenvalues and zeroes of $D$.}

\proclaim{Proposition \thmlbl{evalue}}
For $(\tilde \xi,\lambda)$ within the region of consistent splitting
$\Omega \setminus \{(0,0)\}$, equation \eqnref{gfirstorder} admits 
a nontrivial solution $W\in L^2(x)$ if and only if
$D^\epsilon(\tilde \xi,\lambda)=0$.
\endproclaim

\medskip
{\bf Proof.}
Evidently, the first $K$ columns of the matrix on the 
righthand side of \eqnref{evans} are a basis for the stable 
manifold of \eqnref{gfirstorder} at $x\to +\infty$,
while the final $N-K$ columns are a basis
for the unstable manifold at $x\to -\infty$.
Thus, its determinant vanishes if and only if these
manifolds have nontrivial intersection, and the result follows.
\qed

\medskip
{\bf Convergence of approximate flows.}
We now turn to the crucial question:
Under what circumstances does the Evans function of \eqnref{gfirstorder}
converge as $\epsilon \to 0^+$ to the Evans function of the limiting  
equations at $\epsilon=0$?
The following proposition gives a simple and sharp answer.

\proclaim{Proposition \thmlbl{evanslimit}}
Suppose that, in an $\Omega$-neighborhood of some $(\tilde \xi,\lambda)$:

(i) As $\epsilon \to 0^+$, the asymptotic subspaces 
$S^\epsilon_+$, $U^\epsilon_-$
converge uniformly in angle to $S^0_+$, $U^0_-$, with rate 
$\delta(\epsilon)$: equivalently, for $0<\epsilon \le \epsilon_0$,
their spanning bases satisfy
$$
|v^\epsilon_{j\pm}- v^0_{j\pm}| \le \eta(\epsilon).
\eqnlbl{Sest}
$$

(ii) The coefficient matrices $\BbbA^\epsilon$
converge uniformly exponentially to their limits,
in the sense that, for $0<\epsilon\le \epsilon_0$, 
$$
|(\BbbA^\epsilon - \BbbA^\epsilon_\pm)-
(\BbbA^0- \BbbA^0_\pm)|
\le C_1\eta(\epsilon)e^{-|x|/C_2}.
\eqnlbl{residualest}
$$

Then, on some (possibly smaller) $\Omega$-neighborhood of
$(\tilde \xi,\lambda)$,  the local Evans function
$D^\epsilon $ defined as in Definition \thmref{evans} converges
uniformly as $\epsilon \to 0^+$ to $D^0$, with
$$
|D^\epsilon- D^0|\le
C C_1 C_2 \eta(\epsilon).
\eqnlbl{evansest}
$$
(More precisely, the columns of the defining matrix on the righthand
side of \eqnref{evans} converge uniformly, with the same rate
$CC_1C_2\eta(\epsilon)$.)
\endproclaim

\medskip
{\bf Proof.}
Clearly, it is sufficient to prove the stronger parenthetical assertion
of convergence of individual columns, and for this we may restrict
without loss of generality to the $+$ columns, and the right
half-line $x\ge 0$.
Applying the conjugating transformation $W\to (P^0_+)^{-1}W$ 
for the $\epsilon=0$ equations, we may reduce to the case that
$\BbbA^0$ is constant, and $P^0_-\equiv 0$.
In this case, \eqnref{residualest} becomes just
$$
|(\BbbA^\epsilon - \BbbA^\epsilon_\pm)|
\le C_1\eta(\epsilon)e^{-|x|/C_2},
$$
and we obtain directly from the conjugation lemma, Lemma \thmref{conjugation},
the estimate
$$
|P^\epsilon_+ - P^0_+|=
|P^\epsilon_+ - I|\le CC_1 C_2 \eta(\epsilon) e^{-\theta|x|/C_2},
$$
and in particular
$$
|P^\epsilon_+ - P^0_+|_{|x=0}\le CC_1 C_2 \eta(\epsilon). 
\eqnlbl{zeroPest}
$$
The result now follows, from
\eqnref{evans}, \eqnref{Sest}, and \eqnref{zeroPest}.
\qed
\medskip

{\bf Remark.}  The same argument shows that the proposition is
sharp, since easy examples show that the conjugation lemma itself is sharp.
\medskip

{\bf 2.2. The tracking lemma and reduction of the flow.} 
Next, consider complementary situation of a family of equations of
form \eqnref{gfirstorder} on an $\epsilon$-neighborhood
for which the coefficient
$\BbbA^\epsilon$ does not exhibit uniform exponential
decay to its asymptotic limits, but instead is {\it slowly varying}.
This occurs quite generally for rescaled eigenvalue equations arising in the
study of the large-frequency regime; see, e.g., [GZ,ZH,MZ.1,Z.3].
In the present context, it arises for the {\it unrescaled equations}
in the small-shock strength limit $\epsilon \to 0$.

In this situation, it frequently occurs that not only $\BbbA^\epsilon$
but also certain of its invariant (group) eigenspaces are slowly varying with
$x$, i.e., there exist matrices 
$$
 L^\epsilon=\pmatrix L^\epsilon_1 \\ L^\epsilon_2\endpmatrix(x),
\quad 
R^\epsilon=\pmatrix R^\epsilon_1 &  R^\epsilon_2\endpmatrix (x)
\eqnlbl{LR}
$$
for which $L^\epsilon R^\epsilon(x)\equiv I$
and $|LR'|=|L'R|$ is small relative to
$$
\BbbM^\epsilon:=
L^\epsilon \BbbA^\epsilon R^\epsilon(x) =
\pmatrix M^\epsilon_1 & 0 \\
0 & M^\epsilon_2 \\
\endpmatrix(x),
\eqnlbl{M}
$$
where ``$'$'' as usual denotes $\partial/\partial x$.
In this case, making the change of coordinates $W^\epsilon=R^\epsilon Z$,
we may reduce \eqnref{gfirstorder} to the approximately block-diagonal
equation
$$
{Z^\epsilon}'= \BbbM^\epsilon Z^\epsilon + \delta^\epsilon
\Theta^\epsilon Z^\epsilon,
\eqnlbl{blockdiag}
$$
where $\BbbM^\epsilon$ is as in \eqnref{M},
$\Theta^\epsilon(x)$ are uniformly bounded matrices, 
and $\delta^\epsilon(x)\le \delta(\epsilon)$ is a (relatively) 
small scalar.
A sometimes crucial improvement may be obtained by
arranging that error $\Theta$ vanish on
diagonal blocks, i.e., that $L_j R_j'\equiv 0$; see [Go.1--2,MZ.1].
Here, we shall make use of this observation only in the 
principal,  $2\times 2$ ``Burgers'' block.

Let us assume that such a procedure has been successfully carried out,
and, moreover, that there exists a {\it uniform spectral gap}, 
in the strong sense that
$$
\R (M_1- M_2) \ge \eta(\epsilon) >0
\, \text{\rm for all } x.
\eqnlbl{gap}
$$
In the ``standard'' case that $\BbbM^\epsilon$ are uniformly bounded
and $\eta(\epsilon)$ may be taken independent of $\epsilon$,  
that \eqnref{gap} may be arranged by a suitable coordinate transformation 
provided that there holds the weaker condition
$$
\min \R \sigma(M_1)- \max \R \sigma (M_2) \ge \tilde \eta >0
\, \text{\rm for all } x,
\eqnlbl{weakgap}
$$
for some $\tilde \eta > \eta$.
Then, there holds the following {\it reduction lemma}, a refinement
established in [MZ.1] of the ``tracking lemma'' given in varying
degrees of generality in [GZ,ZH,Z.3].
The new feature of the reduction lemma 
is that it asserts the existence of smooth invariant manifolds
in the vicinity of the first and second block coordinates,
whereas the tracking lemma(s) assert only the existence of
forward and backward attracting cones in the same vicinity
and do not give a reduction in the usual dynamical systems sense.

\proclaim{Proposition \thmlbl{reduction}[MZ.1]}
Consider a system \eqnref{blockdiag} under the gap assumption \eqnref{gap},
with $\Theta^\epsilon$ uniformly bounded for $\epsilon$ sufficiently small.
If, for $0<\epsilon<\epsilon_0$,
the ratio $\delta(\epsilon)/\eta(\epsilon)$ of formal error
vs. gap is sufficiently small relative to the bounds on $\Theta$,
in particular if $\delta(\epsilon)/\eta(\epsilon)\to 0$, then, for 
all $0<\epsilon\le \epsilon_0$, there exist (unique) linear transformations
$\Phi_1^\epsilon(x,\tilde \xi,\lambda)$ and
$\Phi_2^\epsilon(x,\tilde \xi,\lambda)$,
possessing the same regularity with respect to the various parameters
$\epsilon$, $x$, $\tilde \xi$, $\lambda$ as do coefficients
$\BbbM^\epsilon$ and $\delta(\epsilon)\Theta^\epsilon$,
for which the graphs
$\{(Z_1, \Phi^\epsilon_2 Z_1)\}$ and $\{(\Phi^\epsilon_1(Z_2),Z_2)\}$
are invariant under the flow of \eqnref{blockdiag}, and satisfying
$$
|\Phi^\epsilon_1|, \, |\Phi^\epsilon_1| \le C \delta(\epsilon)/\eta(\epsilon)
\, \text{\rm for all } x.
$$
\endproclaim
\medskip
{\bf Proof.}
This may be established by a contraction mapping argument, carried
out for the projectivized ``lifted'' equations governing the
flow of exterior forms;
see Appendix A.3.2.2 (ii) of [MZ.1].
See, e.g.,  [Sat] for a corresponding argument in the case that $M_1$, $M_2$
are scalar;
the lifting to exterior forms essentially reduces the problem to this case.
\qed
\medskip

From Proposition \thmref{reduction}, we obtain reduced flows
$$
{Z_1^\epsilon}'= M_1^\epsilon Z_1^\epsilon + 
\delta^\epsilon \Theta_{12}^\epsilon \Phi_2^\epsilon Z_1^\epsilon
= M_1^\epsilon Z_1^\epsilon + 
\Cal{O}(\delta^\epsilon(x))Z_1^\epsilon
\eqnlbl{reduced1}
$$
and
$$
{Z_2^\epsilon}'= M_2^\epsilon Z_2^\epsilon + 
\delta^\epsilon \Theta_{21}^\epsilon \Phi_1^\epsilon Z_2^\epsilon
= M_2^\epsilon Z_2^\epsilon + 
\Cal{O}(\delta^\epsilon(x))Z_2^\epsilon
\eqnlbl{reduced2}
$$
on the two invariant manifolds described.
Should we arrange that 
$\Theta_{11} \equiv0$ or $\Theta_{22}\equiv 0$, then the error terms
would become 
$\Cal{O}(\delta^\epsilon(x)\delta(\epsilon)/\eta(\epsilon))Z_1^\epsilon$
or 
$\Cal{O}(\delta^\epsilon(x)\delta(\epsilon)/\eta(\epsilon))Z_2^\epsilon$,
respectively; that is, enforcing vanishing of $\Theta$ on a diagonal block
reduces the error term by factor $\delta/\eta$.

\medskip
{\bf Remark \thmlbl{relaxedgap}.}
As pointed out in [MZ.1], the gap condition \eqnref{gap} may
be allowed to fail by order $\alpha$ on an interval of length
order $1/\alpha$, for any fixed $\alpha$.
For the applications of this paper and of [PZ], as for those of [MZ.1],
this means that the gap condition
need only be checked at $x=\pm \infty$.
\medskip
{\bf Remark \thmlbl{errorexp}.}
As described in [MZ.1], provided that $M_1$ and $M_2$ are (uniformly)
bounded and spectrally separated, one can in fact 
repeat the procedure used to construct
\eqnref{blockdiag} to obtain an approximate block-diagonalization
to arbitrarily high order: i.e., an asymptotic expansion in powers of
$\delta(\epsilon)$.
In general, e.g., in the large-amplitude situation considered
in [MZ.1], these higher order terms are not explicitly computable.
However, we point out that, in the present small-amplitude situation,
they may in principle be determined exactly to any desired order, 
using the error expansion for the traveling wave profile to approximate $L'$;
indeed, this can be recognized as just a variation of the usual center 
manifold calculations used to approximate center manifold flow.
In the calculations of this paper and of [PZ], 
we shall not need to perform any such
higher-order corrections: a single iteration will be enough.

\medskip
{\bf 2.3. Reduction to small frequency.}
As a first application of Proposition \thmref{reduction}, we
obtain the following preliminary result analogous to
that of [GJ], reducing the small-amplitude
stability problem to the small frequency, or ``diffusive,'' regime
$|\lambda|<<1$.

\proclaim{Proposition \thmlbl{smallfreq}}
In each of the contexts considered in this paper (one-dimensional 
viscous and relaxation systems, with hypotheses
given below), the nonstable eigenvalues $\R \lambda\ge 0$ 
of linearized operator $L$ 
of \eqnref{ablinear} are restricted to 
$|\lambda|\le r(\epsilon)$, 
$r(\epsilon)\to 0$ as $\epsilon \to 0$, where $\epsilon:=|u_+-u_-|$
denotes shock strength.
\endproclaim

\medskip
{\bf Proof.}
Equivalently, we show, for arbitrary $r>0$, that $|\lambda|\le r$
for $\epsilon $ sufficiently small.
For $|\lambda|>r$ and $\BbbA^\epsilon$ as in 
\eqnref{gfirstorder} denoting the coefficient of the generalized
eigenvalue problem \eqnref{ablinear} written as a first-order
system, there is a uniform spectral gap between
stable and unstable subpaces of $\BbbA^\epsilon(x)$, independent of 
$0<\epsilon\le \epsilon_0$ and $x$; as discussed in [Z.3], this
is a straightforward consequence of the assumed dissipativity 
of the systems (in the sense of Kawashima--Zeng [Kaw,Ze]).
By standard matrix perturbation theory [Kat] plus boundedness (compactness) 
of $\BbbA^\epsilon$, we may therefore deduce the existence of well-conditioned
bases \eqnref{LR} for these subspaces varying smoothly with $\BbbA^\epsilon$,
and therefore satisfying 
$$
|LR'|=\Cal{O}(|(\bar u^\epsilon)'|)=\Cal{O}(\epsilon^2)
$$
by Proposition \thmref{limprofile}.  Moreover, by a proper choice of
basis, we may arrange \eqnref{gap} with $\eta>0$ uniformly bounded,
independent of $0<\epsilon\le \epsilon_0$ and $x$; for details, see 
Appendix A4, [Z.3].

Applying Proposition \thmref{reduction}, we obtain a pair of reduced
systems \eqnref{reduced1} and \eqnref{reduced2} for which the
coefficent matrices (with $\Cal{O}(\epsilon^2)$ error term taken
into account) are, respectively, uniformly positive and negative definite,
and solutions therefore uniformly exponentially grow and decay,
provided that $\epsilon$ is sufficiently small.
It follows that the equations have no solutions bounded at plus
and minus spatial infinity save for the trivial solution $W\equiv 0$.

For $|\lambda|\ge R$, $R>0$ sufficiently
large, on the other hand, the nonexistence of eigenvalues may be established
independent of the size of $\epsilon$
by a rescaling followed by a similar reduction argument;
see [ZH,MZ.1--2]. 
for a multidimensional version, see Lemma 4.38, [Z.3].
\qed

\medskip
{\bf Remark \thmlbl{md}.}
The same argument applies word for word in the multidimensional case,
to yield $|(\tilde \xi, \lambda)|<<1$, where $\tilde \xi$ denotes
the Fourier transform frequency in directions parallel to the shock
front, and $\lambda$ as here denotes the Laplace
transform frequency with respect to time.

\bigskip
\newsection {\bf Section 3. One-dimensional stability for parabolic systems.}
\sectionnumber=3 
\theoremnumber=0
\equationnumber=0
\smallskip
\TagsOnLeft

With these preparations, we are ready to carry out the stability analysis.
We begin with a treatment of small-amplitude stationary profiles
$\bar u^\epsilon(x)$ of one-dimensional strictly
parabolic viscous conservation laws
$$
u_t + f(u)_{x}=(B(u) u_{x})_{x} 
\eqnlbl{1dviscous}
$$
in a neighborhood $\Cal{U}$ of a particular state $u_0$. 
We make the assumptions:
\medskip

(H0)  $f, B\in C^2$ (regularity).
\smallskip
(H1)  $\R \sigma (B(\bU) )>0$ (strict parabolicity).
\smallskip
(H2)  $\sigma(Df)$ real, simple
(strict hyperbolicity).
\smallskip
(H3)  $\R \sigma(-i Df\xi-B \xi^2)
\le -\theta|\xi|^2$, $\xi \in \BbbR$,
for some $\theta>0$  
(linearized stability of constant states).
\smallskip
(H4) $a_p=0$ is a simple eigenvalue of $Df(u_0)$
with left and right eigenvectors $l_p$ and $r_p$,
and $l_p^t D^2f(r_p,r_p)\ne 0$ (genuine nonlinearity of the principal 
characteristic field).
\medskip
Here, $\epsilon>0$ denotes shock strength $|u^\epsilon_+-u^\epsilon_-|$, 
and profiles $\bar u^\epsilon(\cdot)$ are assumed to 
converge as $\epsilon\to 0$ to $u_0$.
Under these assumptions, the center-manifold argument of Majda\&Pego [MP]
verifies the assertion of Proposition \thmref{limprofile}, yielding
convergence after rescaling of $\bar u^\epsilon$ to the standard
Burgers profile \eqnref{burgersprofile}.
In the rest of this section, we shall establish our first main result:

\proclaim{Theorem \thmlbl{pstab}}
Under assumptions (H0)--(H4), profiles $\bar u^\epsilon$
are strongly spectrally stable (and therefore linearly
and nonlinearly orbitally stable [ZH,Z.2])
for $\epsilon$ sufficiently small.
\endproclaim

More precisely, we shall establish the following 
proposition, from which, together  
with Proposition \thmref{modelstab},
Theorem \thmref{pstab} follows as an immediate corollary.

\proclaim{Proposition \thmlbl{equiv}}
Under assumptions (H0)--(H4), profiles $\bar u^\epsilon$
are strongly spectrally stable (and therefore linearly
and nonlinearly orbitally stable [ZH,Z.2])
for $\epsilon$ sufficiently small, if the
standard Burgers profile \eqnref{burgersprofile} is
strongly spectrally stable for $\epsilon=1$ and
only if the Burgers profile is not strongly spectrally
unstable for $\epsilon=1$ in the sense that there
exists an eigenvalue $\R \lambda>0$.
\endproclaim

\medskip
{\bf Remark \thmlbl{vacuous}.}
Of course, the ``only if'' part of Proposition \thmref{equiv}
follows vacuously in the case at hand; the point is that we
obtain this information by our method of proof, independent
of any knowledge of behavior of the the limiting system.
This distinction may be important in future
applications to more general situations.
\medskip
{\bf Remark \thmlbl{extensions}.}
Theorem \thmref{pstab} recovers the spectral stability result
obtained by energy methods in [Go.1--2], and slightly extends it
from the class of $A$, $B$ such that $LBR>0$ for some diagonalizing
bases $L$, $R=L^{-1}$ of $A$ to the class of all ``stable'' pairs
$A$, $B$ in the sense of Majda\&Pego [MP].
(In the case that $A$, $B$
are simultaneously symmetrizable, these two classes are equivalent;
see [MP].)
It is straightforward with some additional bookkeeping
(see [Z.3,HuZ]) 
to extend our results to the case, 
to which the nonlinear analyses of [ZH,Z.2] also apply,
that (H2) is relaxed to: 
\medskip
(H2') $\sigma (Df)$ real, semisimple.

\medskip

{\bf 3.1. The case $B\equiv I$.}
To illustrate the method, we first carry out the proof of
Proposition \thmref{equiv} and Theorem \thmref{pstab} 
in the setting of identity viscosity $B\equiv I$, 
for which the associated linear algebra is particularly simple.

Linearizing about $\bar u^\epsilon$, we obtain the
family of eigenvalue equations
$$
\tilde w''=(A^\epsilon \tilde w)'+\lambda \tilde w,
\eqnlbl{Iunint}
$$
where
$$
A^\epsilon(x):= Df(\bar u^\epsilon(x)).
\eqnlbl{IA}
$$
Following Goodman [Go.1], we consider instead
the more favorable ``integrated equations''
$$
w''=A^\epsilon w'+\lambda w
\eqnlbl{Iint}
$$
satisfied by the anti-derivative $w(x):=\int_{-\infty}^{x}\tilde w$ 
of $\tilde w$.
For $\lambda \in \{\R \lambda \ge 0\}\setminus \{0\}$,
it is easily seen that
\eqnref{Iunint} possesses $L^2$ solutions if and only
if does \eqnref{Iint} [ZH]: indeed, the $L^2$ solutions of \eqnref{Iunint}
consist exactly of the derivatives of those of \eqnref{Iint}.
Unlike \eqnref{Iunint}, which has a single zero eigenfunction 
$\tilde w=({\bar u}^\epsilon)'$, \eqnref{Iint} possesses no zero 
($L^2$, i.e., decaying) eigenfunctions.
For, the only candidate, $w=\bar u^\epsilon - u_-^\epsilon$, does
not decay as $x\to +\infty$, but converges to $u^\epsilon_+-u^\epsilon_-$.

\medskip
{\bf Remark \thmlbl{flux}.}
It is worth noting that equations \eqnref{Iint} may be obtained by
another route.  As pointed out by Goodman [Go.2],
the flux variable $w:=\tilde w'- A \tilde w$ satisfies the same equation.
\medskip

Writing \eqnref{Iint} as the family of first-order systems
$$
W'=\BbbA^\epsilon(x,\lambda)W,
\eqnlbl{1dfirstorder}
$$
where
$$
W:=\pmatrix w \\ w'\endpmatrix, \quad
\BbbA^\epsilon :=
\pmatrix
0 & I\\
\lambda & A^\epsilon\\
\endpmatrix,
\eqnlbl{Icoeff}
$$
we are ready to begin our analysis.

\proclaim{Lemma \thmlbl{Ivectors}}
The eigenvalues $\mu$ and associated left and right eigenvectors 
$L$, $R$ of $\BbbA^\epsilon(x)$ in \eqnref{Icoeff} may be denoted as
$$
\mu_j^\pm(x)= \frac{
a_j \pm \sqrt{ a_j^2 + 4 \lambda} }
{2},
\eqnlbl{mu}
$$
$$
L_j^\pm(x)=\Big( \lambda l_j/(\lambda + {\mu_j^\pm}^2), 
\mu_j^\pm l_j/(\lambda+ {\mu_j^\pm}^2 ) \Big),
\quad
R_j^\pm(x)=\pmatrix r_j \\ \mu_j^\pm r_j \endpmatrix,
\eqnlbl{LR}
$$
where $a_j(x)$ denote the eigenvalues of $A^\epsilon(x)$ and $l(x)$, $r(x)$
associated left and right eigenvectors, normalized so that 
$l_j r_k= \delta^j_k$
(here, we have suppressed
the $\epsilon$ dependence for the sake of readability).
\endproclaim
\medskip

By (H2) and (H4), we may order the (real) eigenvalues of $A^\epsilon$
as 
$$
a_1, \dots, a_{p-1}\le -\theta<0,
\eqnlbl{ntrans}
$$
$$
a_p \sim \epsilon \bar \eta(\Lambda \epsilon x/\beta
= \Cal{O}(\epsilon),
\eqnlbl{principal}
$$
and
$$
0< \theta\le a_{p+1}, \dots, a_{n},
\eqnlbl{ntrans}
$$
where $\Lambda$ and $\beta$ denote the genuine nonlinearity and
effective diffusion coefficients described in the introduction,
and $\bar \eta$ the standard Burgers profile \eqnref{burgers}.
Like the principal eigenvalue $a_p$, the transverse eigenvalues
$a_j$, $j\ne p$ vary within an $\Cal{O}(\epsilon)$ neighborhood
of the corresponding eigenvalues of $Df(u_0)$.

Evidently, for the indices $j\ne p$ such that $a_j$ is bounded
from zero, the eigenvalues $\mu_j^\pm$ are jointly analytic in
$a_j$ and thus $A^\epsilon$ and $\lambda$ within neighborhoods of 
$Df(u_0)$, $0$;  a brief calculation shows that $L_j^\pm$
and $R_j^\pm$ are bounded, hence analytic in these parameters as well.
The two-dimensional group eigenspace spanned by $L_p^\pm$
and $R_p^\pm$ likewise varies analytically {\it as subspaces},
as $\BbbA^\epsilon$ and $\lambda$ are varied; however, the
individual eigenvectors that make it up do not (indeed, blow
up for $a_p=0$ and $\lambda\to 0$).
We therefore do not attempt to resolve these subspaces, but just
choose convenient analytically varying spanning bases
$$
L_\eta:= \pmatrix l_p& 0\endpmatrix,
L_z:= \pmatrix 0& l_p\endpmatrix,
\eqnlbl{Lp}
$$
and
$$
R_\eta:= \pmatrix r_p\\ 0\endpmatrix,
R_z:= \pmatrix 0\\ r_p\endpmatrix.
\eqnlbl{Rp}
$$

Summing up, we have:

\proclaim{Lemma \thmlbl{Ibases}}
In the transverse fields $j\ne p$, there hold expansions
$$
\aligned
\mu_j^+&= -\lambda/a_j + \lambda^2 /a_j^3 + \dots\\
\mu_j^-&= a_j + \dots,\\
\endaligned
\eqnlbl{j<p}
$$
for $j<p$ and
$$
\aligned
\mu_j^+&= a_j + \dots,\\
\mu_j^-&= -\lambda/a_j + \lambda^2 /a_j^3 + \dots\\
\endaligned
\eqnlbl{j>p}
$$
for $j>p$, analytic in $\lambda$ and smooth (indeed, analytic) in $a_j$, or, 
through $a_j$, in $A^\epsilon$, on sufficiently small neighborhoods 
of $\lambda=0$ and $a_j=a_j(u_0)$ or $A^\epsilon=Df(u_0)$.
Likewise, the eigenvectors $L_j^\pm$, $R_j^\pm$, $j\ne p$ possess the
same regularity, as do bases $L_\eta$, $L_z$ and $R_\eta$, $R_z$
(defined in \eqnref{Lp}--\eqnref{Rp}) for their complementary two-dimensional
invariant subspaces.
\endproclaim

{\bf Proof.} Taylor expansion/direct calculation.
\qed
\medskip


{\bf Approximate block-diagonalization.}
Following the procedure outlined in Section 2.2, define now 
$$
Z=
\pmatrix \fastvar_- \\
\rho_-\\
\eta\\
z\\
\rho_+\\
\fastvar+\\
\endpmatrix
:= \BbbL^\epsilon W,
\eqnlbl{blockvar}
$$
$$
\BbbL^\epsilon:=\pmatrix
L_{\fastvar_-}\\
L_{\rho_-}\\
L_0\\
L_{\rho_+}\\
L_{\fastvar_+}
\
\endpmatrix,
\quad
\BbbR^\epsilon:=\Big(R_{\fastvar_-},R_{\rho_-},R_0,R_{\rho_+},R_{\fastvar_+}
\Big),
\eqnlbl{ILR}
$$
where
$$
L_{\fastvar_-}:=
\pmatrix 
L_1^+\\
\vdots\\
L_{p-1}^+\\
\endpmatrix,
\quad
R_{\fastvar_-}:=
\Big(R_1^+,
\dots,
R_{p-1}^+,
\Big),
\eqnlbl{ILRfastvar-}
$$
$$
L_{\rho_-}:=
\pmatrix
L_{p+1}^+\\
\vdots\\
L_n^+ \\
\endpmatrix,
\quad
R_{\rho_-}:=
\Big(
R_{p+1}^+,
\dots,
R_n^+ ,
\Big),
\eqnlbl{ILRrho-}
$$
$$
L_0:=
\pmatrix
L_\eta\\
L_z\\
\endpmatrix,
\quad
R_0:=
\Big(
R_\eta,
R_z,
\Big)
\eqnlbl{ILRo}
$$
$$
L_{\rho_+}:=
\pmatrix
L_1^-\\
\vdots\\
L_{p-1}^-\\
\endpmatrix,
\quad
R_{\rho_+} :=
\Big(
R_1^-,
\dots,
R_{p-1}^-,
\Big),
\eqnlbl{ILrho+}
$$
and
$$
L_{\fastvar_+}:=
\pmatrix
L_{p+1}^-\\
\vdots\\
L_n^- \\
\endpmatrix,
\quad
R_{\fastvar_+}:=
\Big(
R_{p+1}^-,
\dots,
R_n^- \Big),
\eqnlbl{ILfastvar+}
$$
$\BbbL^\epsilon \BbbR^\epsilon \equiv I$,
where all coefficients depend on $x$ through $A^\epsilon(\bar u^\epsilon(x))$.

This is just such an approximately block-diagonalizing
transformation as described in Section 2.2:
in particular, because $\BbbL^\epsilon$, $\BbbR^\epsilon$ depend
smoothly on $\BbbA^\epsilon$, we have that
$$
\Theta^\epsilon := |\eta'|^{-1}\BbbL^\epsilon (\BbbR^\epsilon)' =\Cal{O}(1),
\eqnlbl{thetabd}
$$
where 
$\eta:=l_p(u_-)\cdot (\bar u^\epsilon- (u_- + u_+)/2)$ 
as in the introduction.
Thus, in $Z$ coordinates, \eqnref{1dfirstorder}
reduces to the form \eqnref{blockdiag}, with
$$
\aligned
\delta^\epsilon(x)&:= |\eta'|= \Cal{O}(\epsilon^2 e^{-\theta|x|}),
\quad \theta>0, \\
\delta(\epsilon)&=\Cal{O}(\epsilon^2),
\endaligned
\eqnlbl{Idelta}
$$
$$
|\Theta^\epsilon|=\Cal{O}(1), 
\eqnlbl{ITheta}
$$
and
$$
\BbbM^\epsilon:=\diag(
M_{\fastvar_-},
M_{\rho-},
M_{0},
M_{\rho_+},
M_{\fastvar_+}),
\eqnlbl{IM}
$$
where, for $\R \lambda \ge 0$:
$$
M_{\fastvar_-}:=\diag(\mu_1^+,\dots,\mu_{p-1}^+)\le -\eta <0,
\eqnlbl{Mfastvar-}
$$
$$
M_{\fastvar_+}:=\diag(\mu_{p+1}^-,\dots,\mu_{n}^-)\ge \eta >0,
\eqnlbl{Mfastvar+}
$$
$$
M_{\rho-}:=\diag(\mu_{p+1}^+,\dots,\mu_{n}^+)\le 
-\eta(\R \lambda + |\Im \lambda|^2) \le 0,
\eqnlbl{Mfastvar+}
$$
$$
M_{\rho+}:=\diag(\mu_1^-,\dots,\mu_{p-1}^-)\ge 
\eta(\R \lambda + |\Im \lambda|^2) \ge 0,
\eqnlbl{Mfastvar-}
$$
with also
$$
|M_{\rho_\pm}|=\Cal{O}(|\lambda|),
\eqnlbl{Mmodulus}
$$
and
$$
M_0:=
\pmatrix
0 & 1\\
\lambda & a_p\\
\endpmatrix.
\eqnlbl{M0}
$$

We make a final improvement, following [Go.1--2,MZ.1], by making
the rescaling 
$l_j(x)\to \alpha^{-1}(x)l_j(x)$, $r_j(x)\to r_j(x) \alpha(x) $,
with $\alpha$ satisfying ODE
$$
\alpha'=-l_j r_j' \alpha.
\eqnlbl{alphanorm}
$$
This achieves the normalization
$$
l_j r_j'\equiv 0
\eqnlbl{normal}
$$
while preserving the previously stated properties 
\eqnref{thetabd}--\eqnref{M0}.
As a consequence, we obtain the key property that
$L_0 R_0'\equiv 0$, i.e., error
$\Theta^\epsilon$ vanishes in the $\eta$, $z$ block:
$$
\Theta_{(\eta,z),(\eta,z)}\equiv 0.
\eqnlbl{thetanorm}
$$
We could carry out a similar normalization at the $L_j$, $R_j$
level (with $\alpha_j$ now matrix-valued)
to annihilate all diagonal blocks $\Theta_{jj}$;
see Lemma 4.9 of [MZ.1].
However, we do not require this for our argument.
(In fact, we only require vanishing in the lower
lefthand corner $\Theta_{z,\eta}=0$ of the $\eta$, $z$
block, which is automatic; however, we
wish to point out early on this more general procedure.)

{\bf First reduction.}  Observing that coefficients $M_{\fastvar_-}$
and $M_{\fastvar_+}$ of the ``fast transverse modes''
$\fastvar_-$ and $\fastvar_+$ are separated by a uniform spectral
gap both from each other, and from the ``slow'' modes 
$(\rho_-,\eta,z,\rho_+)^t$, we may 
apply Proposition \thmref{reduction} twice to reduce \eqnref{blockdiag}
to decoupled flows on three invariant manifolds, associated respectively
with $\fastvar_-$,  $\fastvar_+$, and
$(\rho_-,\eta,z,\rho_+)^t$, appearing as
$$
\fastvar_-'= M_{\fastvar_-}\fastvar_- + \Cal{O}(\delta(\epsilon))\fastvar_-,
\eqnlbl{redfastvar-}
$$
$$
\fastvar_+'= M_{\fastvar_+}\fastvar_+ + \Cal{O}(\delta(\epsilon))\fastvar_+,
\eqnlbl{redfastvar+}
$$
and
$$
\pmatrix \rho_-\\ \eta\\ z\\ \rho_+\\ 
\endpmatrix'
=
\pmatrix
M_{\rho_-}&0&0&0\\
0&0&1&0\\
0&\lambda&a_p&0\\
0&0&0& M_{\rho_+}\\
\endpmatrix
\pmatrix \rho_-\\ \eta\\ z\\ \rho_+\\ 
\endpmatrix
+
\delta^\epsilon \tilde \Theta^\epsilon
\pmatrix \rho_-\\ \eta\\ z\\ \rho_+\\
\endpmatrix,
\eqnlbl{red1}
$$
where, by \eqnref{ITheta},
$$
\tilde \Theta^\epsilon:=
\Theta^\epsilon_{(\rho_-,\eta,z,\rho_+),
(\rho_-,\eta,z,\rho_+)}+
\Cal{O}(\delta(\epsilon))\Theta^\epsilon_{
(\rho_-,\eta,z,\rho_+), (\fastvar_-,\fastvar_+)}
$$
is small on the principal block:
$$
\tilde \Theta^\epsilon_{(\eta,z),(\eta,z)}=
\Cal{O}(\delta(\epsilon))=\Cal{O}(\epsilon^2).
\eqnlbl{ItildeTheta}
$$

Observing that the flow is uniformly exponentially
decreasing for \eqnref{redfastvar-} and increasing for \eqnref{redfastvar+},
we find that the only possible decaying solution are the trivial ones
$\fastvar_\pm\equiv 0$, and so we may discard these ``fast transverse'' 
equations, leaving us with \eqnref{red1}.
Rescaling now by $x\to \Lambda\epsilon x$, 
$\lambda\to \lambda/\Lambda^2\epsilon^2$,
$z\to z/\Lambda\epsilon$, we obtain the block-triangular system
of equations \eqnref{rhoeqn}--\eqnref{superslow} described in the introduction,
modulo coefficient errors
$$
\Phi:=\bar\delta^\epsilon\bar \Theta^\epsilon
+ \epsilon \Cal{E}^\epsilon,
\eqnlbl{errors}
$$
where
$$
\bar\delta^\epsilon(x):=
\epsilon^{-1}\delta^\epsilon(x/\epsilon)
=\CalO(\epsilon e^{-\theta |x|}), \quad \theta>0,
\eqnlbl{tildedelta}
$$
with $\bar \Theta^\epsilon$ uniformly bounded and
$\CalO(\epsilon)$ on the principal,
$(\eta,z)$ block), and
$\Cal{E}^\epsilon\equiv 0$ in the $\rho_\pm$ equations
and for the $(\eta,z)$ equations is
$$
\Cal{E}^\epsilon:=
\pmatrix
0&0\\
0& a_p-\bar \eta
\endpmatrix,
\eqnlbl{Eeps}
$$
giving
$$
\Cal{E}^\epsilon_\pm= \pmatrix
0&0\\
0&  \Cal{O}(1)
\endpmatrix, \quad
|\Cal{E}^\epsilon - \Cal{E}^\epsilon_\pm|=\CalO(e^{-\theta |x|}),
\eqnlbl{Edecay}
$$
by Corollary \thmref{limchar}.
Rearranging, we have the error
bounds asserted in the introduction.
Note that the forcing term
$$
\aligned
N(x,\lambda)&=\pmatrix 0\\ n(x,\lambda)\endpmatrix,\\
n&:= \epsilon^{-2}\delta^\epsilon \Theta_{z, \rho}
=\Cal{O}(e^{-\theta|x|}),\\
\endaligned
\eqnlbl{Ndef}
$$
arises in \eqnref{etaeqn} as a result of the rescaling $z\to z/\epsilon$.
\medskip

{\bf Further reductions/Normal forms.}
Loosely following the introduction, we now examine separately each of the
three regimes (with respect to the rescaled variable $\lambda$)
I. $|\lambda|\le C$, II. $C\le |\lambda| \le C\epsilon^{-1}$:
and III. $C\epsilon^{-1}\le |\lambda| \le C^{-1}\epsilon^{-2}$,
where $C>0$ is a sufficiently large constant to be determined later.
(Recall: In Proposition \thmref{smallfreq} we have already disposed of the case
$|\lambda|\ge C^{-1}\epsilon^{-2}$.)

\medskip
{\it Region I.}
On region I, the rescaled $\rho_\pm$, $(\eta,z)$ equations evidently 
have bounded coefficients satisfying condition (ii) of 
Lemma \thmref{evanslimit}, with $\eta(\epsilon)=\CalO(\epsilon)$.
Moreover, the asymptotic coefficients at plus and minus spatial infinity
are block-diagonal, with $M_+$ and $M_-$ positive and negative
definite, respectively, and $M_0=\bar M_0 + \Cal{O}(\epsilon)$,
where
$$
\bar M_0= \pmatrix
0& 1\\
\lambda & \mp 1
\endpmatrix
$$
at $x=\pm \infty$ have uniform spectral gap between stable and
unstable eigenvalues for all $\R \lambda >0$; it follows by standard
matrix perturbation theory that the stable and unstable subspaces
for the asymptotic coefficient matrices of the  $\epsilon$-approximate system
are within $\Cal{O}(\epsilon)$ of those of the limiting, $\epsilon=0$
system, verifying condition (ii) of 
Lemma \thmref{evanslimit}, with $\eta(\epsilon)=\CalO(\epsilon)$.
Thus, we may conclude from the Lemma that the 
approximate system has a local Evans function $D^\epsilon$ 
that is $\CalO(\epsilon)$ close to the local Evans function (with
some chosen basis) $D^0$ for the limiting system with error terms
dropped.  
As described in the introduction,
the limiting $\rho_\pm$ equations decouple, and are readily seen
to admit only trivial $L^2$ solutions, and the limiting $(\eta,z)$ equation
with $\rho_\pm$ set to zero is just the
Burgers eigenvalue equation.  Since the latter is known to admit 
no decaying solutions for $\lambda=0$, we thus have that the limiting
Evans function $D^0$ has no zeroes on Region I if the Burgers equation
is strongly stable on Region I, and has a zero of strictly positive
real part if the Burgers equation is strongly unstable on Region I
(possesses a positive real part eigenvalue).
By uniform convergence/analyticity in $\lambda$, we may conclude
the same regarding the zeroes of the approximate Evans function $D^\epsilon$,
and thus also the eigenvalues of the original operator about the wave.

\medskip
{\it Region II.}
On Region II, the eigenvalues 
$$
\bar \mu_p^\pm(x)
= \frac{\bar\eta \pm \sqrt{ \bar \eta^2 + 4 \lambda}}{2}(x)
$$
of the $(\eta,z)$ (``Burgers'') block $\bar M_0$ satisfy
$$
\aligned
\bar \mu^-_p &\le -\eta|\lambda|^{1/2}\le -C^{1/2}\eta <0,\\
\bar \mu^+_p &\ge +\eta|\lambda|^{1/2} \ge C^{1/2}\eta >0,\\
\endaligned
\eqnlbl{IIgap}
$$
with also $|\bar \mu_p^\pm|=\CalO(|\lambda|^{1/2})$.
Thus, ``balancing'' the reduced equations by the rescaling
$z\to z/|\lambda|^{1/2}$, we obtain well-conditioned
basis vectors for this block.
At the same time, the forcing matrix $N$ of \eqnref{Nbds}
becomes now $\CalO(|\lambda|^{-1/2})$ (note: we have partially
undone the original rescaling $z \to z/\epsilon$), and may
be viewed as an error term of order $\CalO(C^{-1/2})$:
uniformly small, but not vanishing as $\epsilon \to 0$.
Indeed, reviewing the origins of various error terms, we
find that this is the dominant error term {\it within both
of Regions II and III},
with error terms in the $\rho_\pm$ equations retaining
their originally stated form.

Diagonalizing $\bar M_0$ by a well-conditioned change of basis,
we thus obtain (in rescaled coordinates) an approximately 
block-diagonalized system
$$
\pmatrix \rho_-\\ \eta_-\\ \eta_+\\\rho_+\endpmatrix
=
\diag(M_-,\bar \mu_p^-,\bar \mu_p^+, M_+)
\pmatrix \rho_-\\ \eta_-\\ \eta_+\\\rho_+\endpmatrix
+ \Cal{O}(|\lambda|^{-1/2})
\pmatrix \rho_-\\ \eta_-\\ \eta_+\\\rho_+\endpmatrix,
\eqnlbl{temp}
$$
where $\eta_p^\pm$ denote the now-separated Burgers modes.
Observing that the $\bar \mu_p^\pm$ are separated by
uniform spectral gap $\eta |\lambda|^{1/2}$ both from each other
and from the $\CalO(\epsilon |\lambda|)$ blocks $M_\pm$ on all
of Region II (indeed, on Region III as well), we may therefore
apply Lemma \thmref{reduction} as second time to further reduce 
to decoupled equations on three invariant manifolds, 
associated respectively with $\eta_p^-$, $\eta_p^+$, and $(\rho_-,\rho_+)$,
appearing as
$$
(\eta_p^-)'= (\bar \mu_p^-+ \CalO(C^{-1}) \eta_p^-
\eqnlbl{etap-}
$$
$$
(\eta_p^+)'= (\bar \mu_p^++ \CalO(C^{-1}) \eta_p^+.
\eqnlbl{etap+}
$$
and
$$
\pmatrix \rho_-\\ \rho_+\\ 
\endpmatrix'
=
\pmatrix
M_{-}&0 \\
0& M_+\\
\endpmatrix
\pmatrix \rho_-\\ \rho_+\\ 
\endpmatrix
+
\bar \delta^\epsilon \bar \Theta^\epsilon
\pmatrix \rho_-\\ \rho_+\\ 
\endpmatrix,
\eqnlbl{red2}
$$
where $\bar \delta^\epsilon$ {\it is still $\CalO(\epsilon e^{-\theta|x|})$
as in the original rescaling}: note that we are here using in a crucial
way the fact that the error terms in $\rho_\pm$ coordinates alone are
better than the $\CalO(|\lambda|^{-1/2}$ 
error terms in the $\eta_p^\pm$ blocks, which were the ones (recorded
in the uniform description \eqnref{temp}) that led to the
error/gap ratio $|\lambda|^{-1/2}/ |\lambda|^{1/2}\sim |\lambda|^{-1}$
determining the tracking angle in Proposition \thmref{reduction}.

Noting as usual that $\eta_p^\pm$ flows are uniformly exponentially
growing/decaying, so support only trivial decaying solutions, we
may discard the first two equations, leaving us with the $(\rho_-,\rho_+)$
equation alone.  Noting that coefficients $M_\pm =\CalO(\epsilon \lambda)$
for these``superslow'' modes are bounded on Region II, 
we may then apply Lemma \thmref{evanslimit} to obtain that suitably
chosen local Evans functions for this superslow flow converge uniformly
(as $\CalO(\epsilon)$) to that for the limiting equations 
with error terms omitted.
But, just as in the treatment of Region I, the limiting $\rho_\pm$
equations decouple into uniformly exponentially growing/decaying
flows, which evidently support no nontrivial decaying solutions.
Thus, we may conclude, finally, that there exist no decaying solutions
of the reduced equations, and therefore no eigenvalues for the
original equations, for $\lambda$ within Region II, and $\epsilon $
sufficiently small.

\medskip
{\it Region III.}
The treatment of Region III is similar to but somewhat simpler than
the treatment of Region II.  For, as already pointed out, all
steps in the reduction to $\rho_\pm$ equations \eqnref{red2}
remain valid on this region as well.
Noting that, by \eqnref{superslow}, 
$M_+$ and $M_-$ have uniform spectral gap of order
$\theta C^2\epsilon$ on this region, $\theta>0$ fixed
and $C>0$ arbitrarily large, while
$\bar \delta^\epsilon$ remains uniformly $\CalO(\epsilon)$,
so we may in this case apply Lemma \thmref{reduction} to obtain
the further reduction to decoupled equations on two invariant
manifolds associated with modes $\rho_-$ and $\rho_+$, appearing as
$$
\rho_-'= \Big(M_- + \CalO(C^{-2}\bar \delta^\epsilon \bar 
\Theta^\epsilon_{\rho_-}) 
\Big) \rho_-
$$
and
$$
\rho_+'= \Big(M_+ + \CalO(C^{-2}\bar \delta^\epsilon \bar 
\Theta^\epsilon_{\rho_+}) 
\Big) \rho_+.
$$
But, these, having negative and positive definite coefficients
$$
\aligned
M_\pm + \CalO(C^{-2}\bar \delta^\epsilon \bar 
\Theta^\epsilon_{\rho_\pm}) &\gtrless
\pm\theta C^2\epsilon + \CalO(C^{-2} \epsilon )\\
&\gtrless \pm \theta C^2\epsilon/2,
\endaligned
$$
are uniformly exponentially decaying and growing, respectively,
so support no nontrivial decaying solutions.
We may thus conclude that there exist no eigenvalues in Region III,
for $\epsilon$ sufficiently small.

\medskip
{\bf Conclusions.}
Collecting the results for each regime, we obtain the result
of Proposition \thmref{equiv}, and therefore of Theorem 
\thmref{pstab}.  
This completes the analysis of the case $B\equiv I$.
\qed

\medskip
{\bf 3.2. The general strictly parabolic case.}
%
The general case follows quite similarly as in the case $B\equiv I$,
in terms of overall argument structure.  However, there is an important 
new technical aspect introduced here, in the way that we carry out
the linear algebra underlying the initial reduction to ``slow'' 
coordinates $(\rho_\pm, \eta,z)$.

For general $B$, the linearized eigenvalue equations \eqnref{Iunint}
become
$$
(B^\epsilon \tilde w')'=(A^\epsilon  \tilde w)'+\lambda \tilde w,
\eqnlbl{Iunint}
$$
where
$$
B^\epsilon(x):=B(\bar u^\epsilon(x)),
\eqnlbl{B}
$$
and
$$
A^\epsilon(x) v:= Df(\bar u^\epsilon(x))v
-\Big( DB(\bar u^\epsilon (x))v\Big) (\bar u^\epsilon)'(x),
\eqnlbl{A}
$$
and the associated integrated equations become
$$
w''=(B^\epsilon)^{-1}A^\epsilon w'+\lambda (B^\epsilon)^{-1}w.
\eqnlbl{Iint}
$$
We express these as a family of first-order systems
$$
W'=\BbbA^\epsilon(x,\lambda)W,  
\eqnlbl{1dBfirstorder}
$$
where
$$
\BbbA^\epsilon:=
\pmatrix
0 & I\\
\lambda (B^\epsilon)^{-1} & (B^\epsilon)^{-1}A^\epsilon\\
\endpmatrix.
\eqnlbl{Bcoeff}
$$

Note that $A^\epsilon$ lies
within $\CalO(|(\bar u^\epsilon)'|)$ of $dF(\bar u^\epsilon)$,
hence is likewise strictly hyperbolic, with eigenvalues $a_j$
and eigenvectors $l_j$ and $r_j$ lying 
within $\CalO(|(\bar u^\epsilon)'|)$ of those of $dF(\bar u^\epsilon)$.
For purposes of our argument, therefore, $A^\epsilon$ and
$Df(\bar u^\epsilon)$ are indistinguishable; in particular,
$a_j\le -\eta<0$ for $j<p$, $a_j\ge \eta>0$ for $j>p$,
and the principle eigenvalue $a_p$ after rescaling
approaches to the Burgers profile $\bar \eta$ with the same
rate predicted in Proposition \thmref{limchar}.
(This last is the main point, and holds also in the
more general setting alluded to in Remark \thmref{extensions},
with (H2) replaced by (H2').)

For simplicity of exposition, we make the provisional hypothesis:
$$
\sigma(B^{-1}Df(u_0)) \, \text{\rm distinct}.
\eqnlbl{provisional}
$$
This concerns only the transverse modes $j\ne p$, since 
the argument of [MP] shows that zero is a simple eigenvalue
under assumptions (H1), (H3), and is easily removed at the
expense of further bookkeeping.
The next lemma furnishes the linear algebra necessary to
carry out the reduction to slow coordinates in the general case; 
for similar arguments, see [ZH,Z.3,MZ.1].

\proclaim{Lemma \thmlbl{Bbases}}
Associated with each of the transverse fields $j\ne p$ of $A^\epsilon$, 
there exist eigenvalues $\mu_j^\pm$ and associated left and right eigenvectors
$L_j^\pm$ and $R_j^\pm$ of $\BbbA^\epsilon$, 
analytic in $\lambda$ and $(A^\epsilon,B^\epsilon)$ 
on sufficiently small neighborhoods 
of $\lambda=0$ and $(Df(u_0),B(u_0))$,
with expansions
$$
\mu_j^+= -\lambda/a_j + \lambda^2 \beta_j^2 /a_j^3 + \dots,
\eqnlbl{muj<p+}
$$
$$
\aligned
L_j^+(x)&=
\Big( l_j  ,
-l_jB/a_j \Big) + \dots,
\quad
R_j^\pm(x)=\pmatrix r_j \\ 0 \endpmatrix + \dots,\\
\endaligned
\eqnlbl{LRj<p+}
$$
and
$$
\mu_j^-= \gamma_j + \dots,
\eqnlbl{muj<p-}
$$
$$
\aligned
L_j^-(x)&=
\Big( 0, \tilde s_j/\gamma_j\Big) + \dots,
\quad
R_j^-(x)=\pmatrix s_j \\ \gamma_j s_j \endpmatrix + \dots,\\
\endaligned
\eqnlbl{LRj<p-}
$$
for $j<p$, and expansions
$$
\mu_j^+= \gamma_j + \dots,
\eqnlbl{muj>p+}
$$
$$
\aligned
L_j^+(x)&=
\Big( 0, \tilde s_j/\gamma_j\Big) + \dots,
\quad
R_j^\pm(x)=\pmatrix s_j \\ \gamma_j s_j \endpmatrix + \dots,\\
\endaligned
\eqnlbl{LRj>p+}
$$
and
$$
\mu_j^-= -\lambda/a_j + \lambda^2 \beta_j^2 /a_j^3 + \dots
\eqnlbl{muj>p-}
$$
$$
\aligned
L_j^-(x)&=
\Big( l_j  , -l_jB/a_j \Big) + \dots,
\quad
R_j^\pm(x)=\pmatrix r_j \\ 0 \endpmatrix + \dots,\\
\endaligned
\eqnlbl{LRj>p-}
$$
for $j>p$.
Here, 
$$
\gamma_1,\dots, \gamma_{p-1}\le -\eta<0
<\eta\le \gamma_{p+1},\dots, \gamma_{n}
$$
denote eigenvalues of $(B^\epsilon)^{-1}A^\epsilon$,
with associated left and right eigenvectors $\tilde s_j$, $s_j$,
and $\beta_j:= l_j B^\epsilon r_j$.

Dual to the right- and left-invariant subspaces $\Span_{j\ne p}\{R_j\}$
and $\Span_{j\ne p}\{L_j\}$ are two-dimensional left- and right-invariant
subspaces spanned by vectors $L_\eta$, $L_z$ and $R_\eta$, $R_z$
that are analytic in $\lambda$ and $(A^\epsilon,B^\epsilon)$
on the same neighborhoods 
of $\lambda=0$ and $(Df(u_0),B(u_0))$, 
with expansions
$$
L_\eta=\Big(l_p, l_*  \Big)+ \dots,
\quad
L_z=\Big( 0, \tilde s_p \Big)+ \dots,
\eqnlbl{Letaz}
$$
and
$$
R_\eta=\pmatrix r_p \\ 0 \endpmatrix +\dots,
\quad
R_z=\pmatrix r_* \\ s_p \endpmatrix + \dots,
\eqnlbl{Retaz}
$$
where 
$$
l_j r_*=\gamma_j l_js_p, \quad
l_* s_j=-l_p s_j/\gamma_j, \qquad \text{\rm for }\, j\ne p
\eqnlbl{transconditions}
$$
and
$$
l_p r_*+ l_* s_p=0.
\eqnlbl{orthoconditions}
$$
\endproclaim

{\bf Proof.} 
By standard matrix perturbation theory, we have analytic dependence
in $(A^\epsilon, B^\epsilon)$ and  $\xi$ of the eigenvalues 
$\lambda_j$ and eigenvectors $\tilde l_j$, $\tilde r_j$ 
of the symbol 
$$
P^\epsilon(i\xi):= -i\xi A^\epsilon -\xi^2 B^\epsilon,
\eqnlbl{symbol}
$$
with expansions
$$
\aligned
\lambda_j&= -i\xi a_j -\xi^2 \beta_j + \dots,\\
\tilde l_j&= l_j + \dots,\\
\tilde r_j&= r_j + \dots,\\
\endaligned
\eqnlbl{disp}
$$
about $(Df(u_0),B(u_0)$, $0$, where $\beta_j:= l_j B^\epsilon r_j$.
For each $j\ne p$, $a_j$ is bounded from zero, and so 
\eqnref{disp} may be inverted, by the Analytic Implicit Function Theorem,
to yield analytic series for $\mu:= i\xi$ and therefore $\tilde l_j$,
$\tilde r_j$ in terms of $\lambda$,
$(A^\epsilon, B^\epsilon)$. 
This verifies expansions \eqnref{muj<p+}, \eqnref{muj>p-}
for the ``slow'' transverse modes; the corresponding eigenvector
expansions \eqnref{LRj<p+}, \eqnref{LRj>p-} then follow by the relations
$$
\aligned
L_j&=\Big(\lambda\tilde l_j/(\lambda+ \mu_j^2 \beta_j)\Big), 
\mu_j \tilde l_j B^\epsilon/(\lambda+ \mu_j^2 \beta_j), \\
R_j&=(\tilde r_j, \mu_j \tilde r_j), \\
\endaligned
$$
between $\tilde l_j$, $\tilde r_j$ and $L_j$, $R_j$; for similar
arguments.

Expansions \eqnref{muj<p-}, \eqnref{LRj<p-} and
\eqnref{muj>p+}, \eqnref{LRj>p+} for the ``fast'' transverse
modes follow in more straightforward fashion by \eqnref{provisional}
and the analytic dependence of isolated eigenvalues/eigenvectors
of $\BbbA^\epsilon$ at $\lambda=0$.
Finally, analytic dependence of the complementary two-dimensional invariant
subspaces follows by duality, whence (see [Kat])
we obtain the existence of analytic bases $(L_\eta, L_z)$, 
$(R_\eta, R_z)$.
By inspection, the choices described are valid representatives
at $\lambda=0$.
\qed
\medskip

Now, similarly as in the case $B\equiv I$, we make the definitions
\eqnref{blockvar}--\eqnref{ILfastvar+}, to obtain an approximately
block-diagonal system satisfying \eqnref{thetabd}--\eqnref{Mmodulus},
but with \eqnref{M0} now replaced by
$$
\aligned
M_0&:= 
L_0 \BbbA^\epsilon R_0 \\
&={L_0}_{|\lambda=0} \BbbA^\epsilon {R_0}_{|\lambda=0} 
+({L_0}-{L_0}_{|\lambda=0}) \BbbA^\epsilon_{|\lambda=0} {R_0}_{|\lambda=0} \\
&\quad +{L_0}_{|\lambda=0}\BbbA^\epsilon_{|\lambda=0} (R_0-{R_0}_{|\lambda=0})
+ \CalO(|\lambda|^2)
\\
&=
\pmatrix
\lambda l_*B^{-1}r_p & l_ps_p + \lambda l_*B^{-1}r_* + \gamma_p l_* s_p)\\
\lambda \tilde s_p B^{-1}r_p & \gamma_p+ \lambda \tilde s_p B^{-1}r_*\\
\endpmatrix\\
&\quad
+
\pmatrix
\CalO(|\lambda|) & \CalO(|\lambda|)\\
\CalO(|\lambda||\gamma_p|+ |\lambda|^2) & \CalO(|\lambda|)\\
\endpmatrix,
\endaligned
\eqnlbl{BM0}
$$
where, in the crucial lower lefthand corner of the final (error) matrix,
we have used the special
properties
$$
\aligned
\BbbA^\epsilon_{|\lambda=0} {R_\eta}_{|\lambda=0}&=0, \\
{L_z}_{|\lambda=0}\BbbA^\epsilon_{|\lambda=0}&=
\gamma_p  {L_z}_{|\lambda=0}
\endaligned
\eqnlbl{specialprop}
$$
together with the usual 
$$
|R_0-{R_0}_{|\lambda=0}|, \, |L_0-{L_0}_{|\lambda=0}|=\CalO(|\lambda|)
$$
to obtain the stated bound.
The $\CalO(\cdot)$ terms in \eqnref{BM0}
may be expressed more precisely as
$$
\pmatrix
O_1 \lambda & O_2 \lambda\\
O_3\lambda^2 + O_4\gamma_p \lambda & O_5 \lambda\\
\endpmatrix,
\eqnlbl{coeff}
$$
with coefficients $O_j$ analytic in all parameters.

Finally, we observe as in [MZ.1] 
that we may achieve $L_0 R_0'\equiv 0$
by a renormalization $L_0 \to \alpha^{-1}L_0$, $R_0 \to R_0 \alpha$,
where $\alpha$ is a $2\times 2$ matrix defined by ODE
$$
\alpha'= - L_0 R_0' \alpha, \quad \alpha|_{x=0}=I.
\eqnlbl{alphaODE}
$$
At $\lambda=0$, $L_0 R_0'$ is initially upper-triangular, hence
$\alpha$ is as well.  In other words, the form of $L_0$, $R_0$ are
preserved, and \eqnref{alphaODE} 
amounts to the pair of vector normalizations
$$
l_p r_p'\equiv \tilde s_p s_p'\equiv 0
$$
together with the condition
$$
l_p r_*'+ l_*s_p'\equiv 0
$$
restricting choices of $l_*$, $r_*$.
This does not affect any of the previously stated equations or bounds,
but only accomplishes again that error $\Theta$ vanish on the $\eta$, $z$
block:
$$
\Theta_{(\eta,z),(\eta,z)}\equiv 0.
\eqnlbl{Thetabd}
$$

\proclaim{Lemma \thmlbl{BtoI}}
Matrix $M_0$ defined in \eqnref{BM0} satisfies
$$
M_0 - \tilde M_0=
\pmatrix \CalO(|\lambda|) & \CalO(|\lambda| + \epsilon)\\
\CalO(|\lambda|\epsilon)& \CalO(\lambda + \epsilon^2)\\
\endpmatrix
\eqnlbl{BtoI}
$$
and 
$$
(M_0 - M_0(\pm \infty))-  
(\tilde M_0 - \tilde M_0(\pm \infty))=
\pmatrix \CalO(|\lambda|\epsilon e^{-\theta \epsilon |x|}) 
& \CalO((|\lambda| + 1) \epsilon e^{-\theta \epsilon |x|})\\
\CalO(|\lambda|\epsilon e^{-\theta\epsilon |x|})
& \CalO((\lambda \epsilon + \epsilon^2) e^{-\theta\epsilon|x|})\\
\endpmatrix
\eqnlbl{BtoI}
$$
for $x\gtrless 0$, where
$$
\tilde M_0 :=
\pmatrix 0 & 1\\
\lambda/\beta& a_p/\beta
\endpmatrix,
\eqnlbl{tildeM0}
$$
$\beta$ as in the introduction denoting 
$\beta_p(u_0)=\beta_p + \CalO(\epsilon)$.
\endproclaim

{\bf Proof.}
Taylor expanding about $a_p=0$, we have 
$$
s_p=r_p + O_1(a_p) a_p,
$$ 
$$
\tilde s_p= l_p B/\beta_p+
O_2(a_p) a_p,
$$
(up to scalar multipliers)
and
$$
\gamma_p=a_p/\beta_p +
O_3(a_p) a_p^2,
$$
with $O_j$ analytic.
From these expansions, together with the fact that
$A$, $B$, $l_p$, $l_*$, $r_p$, $r_*$, $\tilde s_p$, $s_p$, and $a_p$,
as well as coefficients $O_j$ in \eqnref{coeff},
all approach their asymptotic values as $x\to \pm \infty$ at
rate proportional to 
$|\bar u^\epsilon- u_\pm|\sim \epsilon e^{-\theta\epsilon |x|}$,
the result readily follows.
\qed
\medskip

From this point on, the argument goes almost exactly as in the case $B\equiv I$.
Carrying out the first reduction eliminating fast transverse modes, 
and rescaling by
$x\to \Lambda\epsilon x/\beta$, 
$\lambda\to \beta \lambda/ \Lambda^2\epsilon^2$,
$z\to \beta z/\Lambda\epsilon$, we obtain again the block-triangular system
\eqnref{rhoeqn}--\eqnref{superslow} described in the introduction,
with the errors \eqnref{Einftybds} there described.
The latter errors are slightly larger than occurred in case $B\equiv I$;
however, they are sufficient that the arguments in each of Regions I--III
carry through exactly as before.  This completes the proof of 
one-dimensional spectral stability in the general strictly parabolic case.
\qed

\medskip
{\bf Remark \thmlbl{jordan}.}
The crucial estimate of the lower lefthand entry in
the final (error) matrix of \eqnref{BM0} is a generalization
of the standard observation, in the case of a one-parameter bifurcation 
in $\lambda$ from 
an isolated {\it semisimple block} $L(0)$, $R(0)$ of matrix $\BbbA(\lambda)$,
$L(0)$ and $R(0$ consisting of genuine eigenvectors with a single
common eigenvalue $\mu(0)$, that
$$
L\BbbA R= L(0)\BbbA R(0)+ \CalO(\lambda^2),
\eqnlbl{semisimple}
$$
which follows from the observation that
$$
(d/d\lambda)L\BbbA(0) R(0)+ 
L(0)\BbbA(0) (d/d\lambda) R=
\mu(0)(d/d\lambda)(LR)\equiv 0.
\eqnlbl{Bcalc}
$$
In the case at hand, which is essentially a two-parameter bifurcation 
in $\lambda$ and $a_p$ from a Jordan block
at $\lambda=0$, $a_p=0$, the computation \eqnref{Bcalc} of course breaks
down in general, but {\it does remain valid} for the single,
lower lefthand entry of the block that corresponds to the
pairing of genuine left and right eigenvectors $L_z$ and $R_\eta$,
respectively, to give the stated result of validity to order
$\CalO(|\lambda|^2 + |\lambda||a_p|)$.  This explains the
success of the calculation.
Indeed, the entire calculation could have been done 
by expansion about $\lambda=0$, $a_p=0$, giving slightly
degraded error terms $\CalO(|\lambda|+ |a_p|)$ in other entries,
since these would in fact have been sufficient for the analysis;
however, we prefer to carry $a_p$ within the computation, since
it is no more difficult and gives sharper estimates.

\bigskip
\newsection {\bf Section 4. One-dimensional stability for relaxation systems.}
\sectionnumber=4 
\theoremnumber=0
\equationnumber=0
\smallskip
\TagsOnLeft

Next, we treat small-amplitude
stationary profiles 
of one-dimensional relaxation systems
$$
 \pmatrix u\\ v \endpmatrix_t
 + \pmatrix \tilde f(u,v)\\ \tilde g(u,v) \endpmatrix_x
 = \pmatrix 0\\ q(u,v)
 \endpmatrix,
 \eqnlbl{1drelaxation}
$$
with
$q(u,v^*(u))\equiv 0$,
$\text{\rm Re }\sigma\big(q_v(u,v^*(u))\big)<0$ as in 
\eqnref{eq}--\eqnref{qv}, and $f(u)=\tilde f(u,v^*(u))$.
We make the assumptions:
\medskip
(H0) \quad $\tilde f$, $\tilde g$, $q \in C^2$ (regularity).

\smallskip
(H1) \quad $\displaystyle{
\sigma\left( (D\tilde f,D\tilde g)^t (u,v)\right)}$ real, semi-simple with
constant multiplicity, and uniformly bounded away from $0$ 
(nonstrict hyperbolicity
of relaxation system, plus weak subcharacteristic condition).

\smallskip

(H2) \quad $\displaystyle{
 \sigma\left( Df(u_\pm)\right)}$ real, simple (strict hyperbolicity
of equilibrium system).
\smallskip

(H3) \quad $\displaystyle{
\R \sigma\left(  i\xi\binom{D\tilde f}{D\tilde g} (u_\pm, v_\pm) +
\binom 0{dq} (u_\pm, v_\pm)\right) \leq \frac{-\theta |\xi|^2}{1+|\xi|^2}}$,
\quad $\theta>0$,
for all $\xi\in \BbbR$ (linearized stability of constant states).
\smallskip
(H4) $a_p=0$ is a simple eigenvalue of $Df(u_0)$
with left and right eigenvectors $l_p$ and $r_p$,
and $l_p^t D^2f(r_p,r_p)\ne 0$ (genuine nonlinearity for the equilibrium
system of the principal characteristic field).
\medskip
As in the previous section, $\epsilon>0$ denotes shock strength 
$|u^\epsilon_+-u^\epsilon_-|$, 
and profiles $(\bar u^\epsilon,\bar v^\epsilon)(\cdot)$ are assumed to 
converge as $\epsilon\to 0$ to $(u_0,v_0)$ (where necessarily 
$v_0=v_*(u_0)$), under
which assumptions the center-manifold argument of [MZ.1]
verifies the assertion of Proposition \thmref{limprofile}, yielding
convergence after rescaling of $\bar u^\epsilon$ to the standard
Burgers profile \eqnref{burgersprofile}; see also related
analyses of [Yo.Z,FZe].
In the rest of this section we establish the following analogs
of Theorem \thmref{pstab} and Proposition \thmref{equiv}.

\proclaim{Theorem \thmlbl{rstab}}
Under assumptions (H0)--(H4), profiles $\bar u^\epsilon$
are strongly spectrally stable (and therefore linearly
and nonlinearly orbitally stable [MZ.1])
for $\epsilon$ sufficiently small.
\endproclaim

\proclaim{Proposition \thmlbl{requiv}}
Under assumptions (H0)--(H4), profiles $\bar u^\epsilon$
are strongly spectrally stable (and therefore linearly
and nonlinearly orbitally stable [MZ.1])
for $\epsilon$ sufficiently small, if the
standard Burgers profile \eqnref{burgersprofile} is
strongly spectrally stable for $\epsilon=1$ and
only if the Burgers profile is not strongly spectrally
unstable for $\epsilon=1$ in the sense that there
exists an eigenvalue $\R \lambda>0$.
\endproclaim

Theorem \thmref{rstab} completes the program of Mascia\&Zumbrun [MZ.1],
resolving the 
open question of
stability of general small-amplitude relaxation profiles.

\medskip
{\bf Remark \thmlbl{realvisc}.}
The ``real,'' or partially parabolic case is closely related to
the relaxation case, and may be treated by quite similar analysis,
following the framework set up in Appendix A2 of [Z.3]; see [Z.6].
In the real viscous case, satisfactory small-amplitude stability
results have already been obtained by energy methods [HuZ].
\medskip

{\bf 4.1. The $2\times 2$ case.}
For clarity, we first carry out the $2\times 2$ case $\tilde f$,
$\tilde g$, $q$, $f\in \BbbR^1$,
which was treated by energy methods in [L.3].
The linearized eigenvalue equations about profile
$(\bar u^\epsilon,\bar v^\epsilon)$ are
$$
\Big(\CalA^\epsilon
\pmatrix u\\v\endpmatrix \Big)'= (\CalQ^\epsilon - \lambda I) 
\pmatrix u\\v\endpmatrix,
\eqnlbl{revalue}
$$
where
$$
\CalA^\epsilon(x):= \pmatrix \tilde f_u& \tilde f_v\\
\tilde g_u& \tilde g_v\\ \endpmatrix
(\bar u^\epsilon,\bar v^\epsilon)(x),
\quad
\CalQ^\epsilon(x):=\pmatrix 0 & 0\\ q_u & q_v\endpmatrix
(\bar u^\epsilon,\bar v^\epsilon)(x).
\eqnlbl{AQ}
$$
Following [Z.3], we make the (nonsingular, by (H1)) change of variables
$$
\pmatrix \tilde f\\\tilde g\endpmatrix:= 
\CalA^\epsilon \pmatrix u\\v\endpmatrix
\eqnlbl{fg}
$$
to obtain a first-order system
$$
\pmatrix \tilde f\\\tilde g\endpmatrix'= 
(\CalQ^\epsilon - \lambda I) (\CalA^\epsilon)^{-1}
\pmatrix \tilde f\\\tilde g\endpmatrix
\eqnlbl{fluxevalue}
$$
in nondivergence form.

So far, this is analogous to the ``flux transform'' 
of [Go.2].  If we were now to convert to a single,
second-order equation for $\tilde f$, then the analogy
would be exact, since the equations would then clearly
not support any $L^2$ (decaying) solutions at $\lambda=0$.
Instead, we take a different approach, keeping
the convenient variables $\tilde f, \tilde g$, and rescaling as
$$
\tilde f\to \tilde f/\lambda,\quad
\tilde g \to \tilde g,
\eqnlbl{fluxtoint}
$$
to obtain 
$$
\pmatrix \tilde f\\\tilde g\endpmatrix'= 
\diag(\lambda^{-1},1)
(\CalQ^\epsilon - \lambda I) (\CalA^\epsilon)^{-1}
\diag(\lambda,1)
\pmatrix \tilde f\\\tilde g\endpmatrix
\eqnlbl{rfluxevalue}
$$
in place of \eqnref{fluxevalue}.
Since $\tilde f'=\lambda u$ in the original equations \eqnref{revalue},
we see that the rescaled $\tilde f$ variable is, in the case of
decaying solutions, given by
$ \tilde f(x)=\int_{-\infty}^x u(z)dz $,
i.e., is exactly the {\it integrated variable} of 
standard energy analyses [L.3,Liu,HuZ]; in particular,
the rescaled equations \eqnref{rfluxevalue} have the
favorable property that they do not support $L^2$
(decaying) solutions at $\lambda=0$.

{\bf Remark \thmlbl{balanced}.}
Equations \eqnref{rfluxevalue} are neither 
pure flux nor pure integrated form, but a new form intermediate
to the two, which we might be called ``balanced flux form.''  
This form will prove useful again in the multidimensional
case treated in [PZ].
(All three forms (flux, integrated, and balanced flux) 
coincide in the one-dimensional viscous case, so this
distinction is not apparent, or necessary, here.)
\medskip

By direct calculation, we have 
$$
\CalA^{-1}=
(\det \CalA)^{-1}
\pmatrix \tilde g_v & -\tilde f_v\\
-\tilde g_u & \tilde f_u
\endpmatrix,
$$
$$
\aligned
\CalQ\CalA^{-1}&=
(\det \CalA)^{-1}
\pmatrix 0 & 0\\
q_u \tilde g_v-q_v \tilde g_u &  -q_u\tilde f_v + q_v \tilde f_u 
\endpmatrix\\
&= \frac{q_v}{\det \CalA}
\pmatrix 0 & 0\\
-g_u &   a_1
\endpmatrix,
\endaligned
$$
where $ g(u):= g(u,v_*(u))$.
Combining, we find that
$$
\aligned
(\CalQ - \lambda)\CalA^{-1}&=
(\det \CalA)^{-1}
\Big(
q_v
\pmatrix 0 & 0\\
-g_u &   a_1
\endpmatrix
- \lambda \pmatrix \tilde g_v & -\tilde f_v\\
-\tilde g_u & \tilde f_u
\endpmatrix
\Big)\\
&=
\frac{q_v }{\det \CalA}
\pmatrix 0 & \lambda \tilde f_v/ q_v\\
-g_u &   a_1
\endpmatrix
+
\pmatrix \CalO(|\lambda|) & 0\\
\CalO(|\lambda|) & \CalO(|\lambda|) 
\endpmatrix,\\
\endaligned
$$
and therefore
$$
\aligned
\diag(\lambda^{-1},1)(\CalQ - \lambda)&\CalA^{-1}\diag(\lambda,1)=\\
&\frac{q_v }{\det \CalA}
\pmatrix 0 & \tilde f_v/ q_v\\
-g_u \lambda &   a_1
\endpmatrix
+
\pmatrix \CalO(|\lambda|) & 0\\
\CalO(|\lambda|^2) & \CalO(|\lambda|) 
\endpmatrix\\
\endaligned
$$
in \eqnref{rfluxevalue}.

Let 
$$
\beta_1 :=
     -\tilde f_v q_v^{-1}(g_u - q_v^{-1}q_u f_u) 
     =-\tilde f_v q_v^{-1}(g_u - q_v^{-1}q_u a_1) 
\eqnlbl{betadef}
$$ 
denote the effective, ``Chapman-Enskog viscosity'' defined 
in \eqnref{Bstar}. 
A consequence of (H3) is that $\beta>0$; see, e.g. [Ze],
or Appendix A.2 of [MZ.1].

\proclaim{Lemma \thmlbl{2x2facts}}
There hold
$
\det\CalA= -\tilde f_v gu + \tilde g_v a_1
= -\tilde f_v g_u + \CalO(\epsilon)
$
and
$$
\aligned
\beta_1= -q_v^{-1}\tilde f_v gu + \tilde f_v q_v^{-2} q_u a_1
&= -q_v^{-1}\det \CalA+(q_v^{-1}\tilde g_v + \tilde f_v q_v^{-2} q_u) a_1\\
&= -q_v^{-1}\det \CalA + \CalO(\epsilon).\\
\endaligned
$$
In particular, $\det \CalA<0$ and
$\sgn \tilde f_v=\sgn g_u \ne 0$,
for $\epsilon$ sufficiently small.
\endproclaim

\medskip
{\bf Proof.}
A column operation reduces $\det \CalA$ to 
$
\det \pmatrix
f_u & \tilde f_v\\
g_u & \tilde g_v\\
\endpmatrix,
$
from which the first assertion follows.
Combining with \eqnref{betadef}, we obtain the second assertion.
\qed

\medskip
{\bf Remark \thmlbl{subchar}.}
Condition $\det \CalA <0 $ for $a_1$ near zero
may be recognized as the subcharacteristic condition, which
in the $2\times 2$ case is equivalent to (H3).
Condition $\tilde f_v\ne 0$ is the ``genuine coupling'' assumption
of Liu [L.3]; as we have just shown, it is in fact a consequence of (H3), 
and not a separate assumption.
\medskip

Now, setting 
$$
\beta:= \beta_1\CalA/q_v(0),
$$
%
rescaling $x\to \Lambda \epsilon x/\beta$, 
$\lambda \to \beta \lambda/\Lambda^2\epsilon^2$,
and balancing $\tilde f \to  -g_u \beta \tilde f/ \Lambda \epsilon $,
we obtain the normal form
$$
\pmatrix \tilde f\\ 
\tilde g\\
\endpmatrix '=
\pmatrix 0 & 1 \\
\lambda & \bar a_1
\endpmatrix 
\pmatrix \tilde f\\ 
\tilde g\\
\endpmatrix 
+
\Phi
\pmatrix \tilde f\\ 
\tilde g\\
\endpmatrix ,
$$
where, by Lemma \thmref{2x2facts}, error term $\Phi$ satisfies bounds 
\eqnref{Phibds}--\eqnref{Einftybds} of the introduction: that is,
just an approximate Burgers eigenvalue equation.
From this form we may immediately conclude the stability results, 
by a simplified (because transverse, $\rho$ terms are absent) version
of the arguments of the previous section.
This completes the analysis of the $2\times 2$ case.
\qed

\medskip

{\bf 4.2. The general case.}
Now, consider the case of general $n$, $r\ge 1$,
starting as in the previous subsection
with the linearized eigenvalue equation written
in the ``balanced flux variables'' introduced in 
the previous subsection:
$$
W'= \BbbA^\epsilon W,
\eqnlbl{grevalue}
$$
$$
\aligned
\BbbA^\epsilon&:=
\diag(\lambda^{-1}I_n, I_r)
(\CalQ^\epsilon - \lambda I) (\CalA^\epsilon)^{-1}
\diag(\lambda I_n, I_r)\\
&= 
\pmatrix I_n & 0\\
q_u & q_v
\endpmatrix
(\CalA^\epsilon)^{-1}
\pmatrix \lambda I_n & 0\\
0 & I_r
\endpmatrix
\endaligned
\eqnlbl{grfluxevalue}
$$
$\CalA^\epsilon$, $\CalQ^\epsilon$ as defined in \eqnref{AQ}.

Define 
$$
(\tilde H, H):=(q_u,q_v) \CalA^{-1},
\eqnlbl{H}
$$
$\tilde H\in \BbbR^{r\times n}$, $H\in \BbbR^{r\times r}$.
As discussed in the treatment of existence given in Appendix
A.1 of [MZ.1], the $r\times r$ matrix $H$
has a special role as the coefficient matrix of the linearized traveling
wave ODE expressed in flux coordinates.
For simplicity of exposition,
similarly as in the viscous case, we make the provisional hypothesis:
\medskip
(A1) \quad 
$\sigma H(u_0)$ distinct;
\medskip
\noindent 
again, this is easily removed at the expense of further bookkeeping.
A consequence of dissipativity, (H3), is that on a sufficiently
small neighborhood of the base point $u_0$, $H(u)$ has a single eigenvalue 
$$
\gamma_{k}=\CalO(\epsilon),
\eqnlbl{gammak}
$$
with the $(r-1)$ remaining
eigenvalues uniformly bounded from zero as $\epsilon \to 0$:
$$
\gamma_{1}, \dots, \gamma_{k-1}\le-\eta<0
< \eta \le \gamma_{k+1}, \dots, \gamma_{r};
\eqnlbl{gammak}
$$
we denote the associated left and right eigenvectors by $\tilde s_j$
and $s_j$, respectively.

At $\lambda=0$, $\BbbA^\epsilon$ reduces to 
$$
\pmatrix 0 & E\\
0 & H
\endpmatrix,
\eqnlbl{zerored}
$$
where
$$
(\tilde E, E) :=
\Big(-I_n,0 \Big) (\CalA^\epsilon)^{-1}.
\eqnlbl{E}
$$
More generally, we have
$$
\BbbA^\epsilon(x,\lambda)=
\pmatrix \lambda \tilde E & E\\
\lambda \tilde H+ \lambda^2 \tilde F & H + \lambda F
\endpmatrix,
\eqnlbl{fullred}
$$
where $\tilde E$, $E$, $\tilde H$, $H$ are as defined in \eqnref{E}
and \eqnref{H}, and
$$
(\tilde F, F) :=
\Big(0,-I_r \Big) (\CalA^\epsilon)^{-1}.
\eqnlbl{rF}
$$

\proclaim{Lemma \thmlbl{Rbases}}
Associated with each eigenvalue $\gamma_j$, $j\ne k$ of $H$, 
there exist ``fast'' 
eigenvalues $\mu_j^f$ and associated left and right eigenvectors
$L_j^f$ and $R_j^f$ of $\BbbA^\epsilon$, 
analytic in $\lambda$ and $(\CalA^\epsilon,\CalQ^\epsilon)$ 
on sufficiently small neighborhoods 
of $\lambda=0$ and $(\CalA^\epsilon,\CalQ^\epsilon)(u_0,v_0)$, 
with expansions
$$
\mu_j^f= \gamma_j + \dots, 
\eqnlbl{rfastmu}
$$
and
$$
L_j^f= \Big(0, \tilde s_j \Big)+\dots,\quad
R_j^f= \pmatrix \gamma_j^{-1}E s_j\\s_j\endpmatrix +\dots,
\eqnlbl{rfastLR}
$$
where $\gamma_j$, $\tilde s_j$, $s_j$ $\tilde s_j$ 
are eigenvalues and associated left and right eigenvectors of $H$,
and $E$ is as defined in \eqnref{E}.

Likewise, associated with each of the transverse fields 
$j\ne p$ of $A^\epsilon$, 
there exist ``slow'' eigenvalues 
$\mu_j^s$ and associated left and right eigenvectors
$L_j^s$ and $R_j^s$ of $\BbbA^\epsilon$, 
analytic in $\lambda$ and $(\CalA^\epsilon,\CalQ^\epsilon)$ 
on sufficiently small neighborhoods 
of $\lambda=0$ and $(\CalA^\epsilon,\CalQ^\epsilon)(u_0,v_0)$, 
of $\lambda=0$ and $(\CalA, \CalQ)(u_0,v_0))$,
with expansions
$$
\mu_j^s= -\lambda/a_j + \lambda^2 \beta_j^2 /a_j^3 + \dots,
\eqnlbl{rslowmu}
$$
$$
\aligned
L_j^s(x)&=
\Big( l_j  , \tilde t_j \Big) + \dots,
\quad
R_j^s(x)=\pmatrix r_j \\ 0 \endpmatrix + \dots,\\
\endaligned
\eqnlbl{rslowLR}
$$
where $\beta_j:= l_j B_* r_j$, with $B_*$ the effective, ``Chapman-Enskog''
viscosity defined in \eqnref{Bstar} of the introduction,
and $\tilde t_j$ satisfies
$$
l_j E + \tilde t_j H=0.
\eqnlbl{tildetrel}
$$

Dual to the right- and left-invariant subspaces 
$\Span_{j\ne k}\{R^f_j\}\oplus \Span_{j\ne p}\{R_j^s\}$
and 
$\Span_{j\ne k}\{L_j^f\}\oplus \Span_{j\ne p}\{L_j^s\}$ 
are two-dimensional left- and right-invariant
subspaces spanned by vectors $L_\eta$, $L_z$ and $R_\eta$, $R_z$
analytic in $\lambda$ and $(\CalA^\epsilon,\CalQ^\epsilon)$ 
on the same neighborhoods 
of $\lambda=0$ and $(\CalA^\epsilon,\CalQ^\epsilon)(u_0,v_0)$,
with expansions
$$
L_\eta=\Big(l_p, l_*  \Big)+ \dots,
\quad
L_z=\Big( 0, \tilde s_p \Big)+ \dots,
\eqnlbl{Letaz}
$$
and
$$
R_\eta=\pmatrix r_p \\ 0 \endpmatrix +\dots,
\quad
R_z=\pmatrix r_* \\ s_p \endpmatrix + \dots,
\eqnlbl{Retaz}
$$
where 
$$
l_j r_*=-\tilde t_j s_p, \quad
l_* s_j=-\gamma_j^{-1}l_p E s_j, \qquad \text{\rm for }\, j\ne p
\eqnlbl{transconditions}
$$
and
$$
l_p r_*+ l_* s_p=0.
\eqnlbl{orthoconditions}
$$
\endproclaim

\medskip
{\bf Proof.}
This follows similarly as in the proof of Lemma \thmref{Bbases},
namely, analytic variation of the fast transverse modes follows by their
spectral separation, while analytic variation of slow transverse modes
follows by inversion (and appropriate conjugation) 
of the associated dispersion relations
$$
\aligned
\lambda_j&= -i\xi a_j -\xi^2 \beta_j + \dots,\\
\tilde \CalL_j&= \CalL_j + \dots,\\
\tilde \CalR_j&= \CalR_j + \dots,\\
\endaligned
\eqnlbl{disp}
$$
of the symbol 
$ P^\epsilon(i\xi):= \CalQ^\epsilon -i\xi \CalA^\epsilon.$
The calculation of expansions \eqnref{disp} may be found, e.g.,
in Appendix A.2 of [Z.3] or Appendix A.2 of [MZ.1].
Analytic variation of the two-dimensional complementary subspaces
again follows by duality.
\qed
\medskip

As in the treatment of the viscous case, we make the definitions
\eqnref{blockvar}--\eqnref{ILfastvar+}, to obtain an approximately
block-diagonal system satisfying \eqnref{thetabd}--\eqnref{Mmodulus},
with \eqnref{M0} replaced by the following analog of \eqnref{BM0}:

$$
\aligned
M_0&:= 
L_0 \BbbA^\epsilon R_0 \\
&={L_0}_{|\lambda=0} \BbbA^\epsilon {R_0}_{|\lambda=0} 
+({L_0}-{L_0}_{|\lambda=0}) \BbbA^\epsilon_{|\lambda=0} {R_0}_{|\lambda=0} \\
&\quad +{L_0}_{|\lambda=0}\BbbA^\epsilon_{|\lambda=0} (R_0-{R_0}_{|\lambda=0})
+ \CalO(|\lambda|^2)
\\
&=
\pmatrix
\CalO(\lambda) & l_pE s_p + \gamma_p l_* s_p + \CalO(\lambda)\\
\lambda \tilde s_p \tilde H r_p & \gamma_p+ \CalO(\lambda)\\
\endpmatrix,\\
&\quad
+
\pmatrix
\CalO(|\lambda|) & \CalO(|\lambda|)\\
\CalO(|\lambda||\gamma_p|+ |\lambda|^2) & \CalO(|\lambda|)\\
\endpmatrix,
\endaligned
\eqnlbl{rM0}
$$
where all $\CalO(\cdot)$ terms are analytic in $\lambda$, $\CalA$, $\CalQ$.
Here, as in the viscous case, we have used \eqnref{specialprop}
to improve the crucial bound in the lower lefthand corner
of the final, error matrix; more precisely, we have description \eqnref{coeff}
for this error bound, with coefficients $O_j$ analytic in all parameters.
As in the viscous case, we may achieve
$L_0 R_0'\equiv 0$ by a renormalization 
leaving the form of the equations otherwise unchanged;
thus, error $\Theta$ vanishes on the $\eta$, $z$ block:
$$
\Theta_{(\eta,z),(\eta,z)}\equiv 0.
\eqnlbl{Thetabd}
$$

\proclaim{Lemma \thmlbl{rBtoI}}
Matrix $M_0$ defined in \eqnref{rM0} satisfies
$$
M_0 - \tilde M_0=
\pmatrix \CalO(|\lambda|) & \CalO(|\lambda| + \epsilon)\\
\CalO(|\lambda|\epsilon)& \CalO(\lambda + \epsilon^2)\\
\endpmatrix
\eqnlbl{BtoI}
$$
and 
$$
(M_0 - M_0(\pm \infty))-  
(\tilde M_0 - \tilde M_0(\pm \infty))=
\pmatrix \CalO(|\lambda|\epsilon e^{-\theta \epsilon |x|}) 
& \CalO((|\lambda| + 1) \epsilon e^{-\theta \epsilon |x|})\\
\CalO(|\lambda|\epsilon e^{-\theta\epsilon |x|})
& \CalO((\lambda \epsilon + \epsilon^2) e^{-\theta\epsilon|x|})\\
\endpmatrix
\eqnlbl{BtoI}
$$
for $x\gtrless 0$, where
$$
\tilde M_0 :=
\pmatrix 0 & 1\\
\lambda/\beta& a_p/\beta
\endpmatrix,
\eqnlbl{tildeM0}
$$
$\beta:=\beta_p(0)$, and as usual $\beta_j:=l_j B_*r_j$ with 
$B_*$ as defined in \eqnref{Bstar}.
\endproclaim

{\bf Proof.}
Taylor expanding about $a_p=0$, we have 
up to scalar multipliers
(see Appendix A.1, [MZ.1], Claim just
below equation (8.7)):
$$
s_p=g_u r_p + O_1(|a_p|) a_p,
\eqnlbl{sexp}
$$ 
$$
\tilde s_p= -l_p \tilde f_v q_v^{-1}/\beta_p+
O_2(|a_p|) a_p,
\eqnlbl{tildesexp}
$$
and (by a standard matrix perturbation calculation
for a simple eigenvector: see (8.14)--(8.15) of
Appendix A.1, [MZ.1]).
$$
\gamma_p=a_p/\beta_p +
O_3(|a_p|) a_p^2,
\eqnlbl{gammaexp}
$$
where $O_j$ are analytic. 

Furthermore, at $a_p=0$, there hold the relations:
$$
\Big(\tilde E, E\Big)\pmatrix  0 \\ s_p \endpmatrix
= r_p,
\eqnlbl{EErel}
$$
$$
\Big(\tilde F, F\Big)\pmatrix  0 \\ s_p \endpmatrix
= -q_v^{-1}q_u r_p
\eqnlbl{FFrel}
$$
(see (8.8), Appendix A.1, [MZ.1]),
and
$$
\tilde s_p \Big(\tilde H, H\Big)
= 
\Big(-l_p/\beta,  0 \Big)
\eqnlbl{HHrel}
$$
(see (8.10), Appendix A.1, [MZ.1]),
from which we obtain the key facts:
$$
Es_p=r_p, \quad
Fs_p=-q_v^{-1}q_u r_p, \quad
\tilde s_p \tilde H=-l_p/\beta.
\eqnlbl{keyfacts}
$$

From expansions \eqnref{sexp}--\eqnref{gammaexp}, 
the relations \eqnref{keyfacts}, and the fact that
$A$, $\beta_p$, $l_p$, $l_*$, $r_p$, $r_*$, $\tilde s_p$, $s_p$, and $a_p$,
as well as coefficients $O_j$ in \eqnref{coeff},
all approach their asymptotic values as $x\to \pm \infty$ at
rate proportional to 
$|\bar u^\epsilon- u_\pm|\sim \epsilon e^{-\theta\epsilon |x|})$,
the result readily follows.
\qed
\medskip
{\bf Remark.} 
Note that the natural choice of basis vectors $L_z$ and $R_z$
given in Lemma \thmref{Rbases} automatically accomplishes the
desired rescaling to the viscous case, so that we need not do
this ``by hand'' as in the treatment of the $2\times 2$ case above.
\medskip

From this point, the argument goes through
exactly as in the general viscous case.
Carrying out the first reduction eliminating the $r-1$
fast transverse modes, and rescaling by
$x\to \Lambda\epsilon x/\beta$, 
$\lambda\to \beta \lambda/ \Lambda^2\epsilon^2$,
$z\to \beta z/\Lambda\epsilon$, we obtain as in the
general viscous case the block-triangular system
\eqnref{rhoeqn}--\eqnref{superslow} described in the introduction,
with errors \eqnref{Einftybds}, and so the remaining arguments
go through as before.
This completes the proof of 
one-dimensional spectral stability in the general relaxation case.
\qed

\bigskip
\newsection {\bf Section 5. Multidimensional stability for parabolic
systems.}
\sectionnumber=5 
\theoremnumber=0
\equationnumber=0
\smallskip
\TagsOnLeft

Finally, we briefly discuss the extension of the above analysis
to the multidimensional viscous case; details will be given in [PZ].
Consider a sequence of planar stationary profiles $\bar u^\epsilon(x_1)$ 
of a strictly parabolic viscous conservation law of form
$$
u_t + \sum_j f^j(u)_{x_j}=\sum_{jk}(B^{jk}(u) u_{x_k})_{x_j},
\eqnlbl{viscous}
$$
lying in a neighborhood $\Cal{U}$ of a particular state $u_0$, where
$x\in \BbbR^d$, $u$, $f^j \in \Bbb{R}^n$, and $B^{jk}\in \BbbR^{n\times n}$.
We make the assumptions:
\medskip

(H0)  $f^j, B^{jk}\in C^2$ (regularity).
\smallskip
(H1)  $\R \sigma (\sum B^{jk}(\bU) \xi_j \xi_k)>0$ for $\xi \in \BbbR^d$
(strict parabolicity).
\smallskip
(H2)  There exists $A^0(\cdot)$, symmetric positive definite
and smoothly depending on $u$,
such that $A^0 Df^j(u)$ is symmetric for all $1\le j\le d$.
(simultaneous symmetrizability, $\Rightarrow$ nonstrict hyperbolicity).
\smallskip
(H3)  $\R \sigma(-i\sum Df^j(u_\pm)\xi_j-\sum B^{jk}(u_\pm)\xi_j
\xi_k)\le -\theta|\xi|^2$, $\xi \in \BbbR^d$,
for some $\theta>0$  
(linearized stability of constant states).
\smallskip
(H4) (i) $a_p=0$ is a simple eigenvalue of $Df^1(u_0)$
with left and right eigenvectors $l_p$ and $r_p$,
and $l_p^t D^2f^1(r_p,r_p)\ne 0$ (genuine nonlinearity of the principal 
characteristic field in the normal spatial direction $x_1$).

(ii) $r_p$ is not an eigenvector of any $Df^{\tilde \xi}:=
\sum_{j\ne 1}\xi_j Df^j$ for $\tilde \xi \in \BbbR^{d-1}\ne 0$,
and, in the intermediate case $1<p<n$:
$$
\langle r_p,
A^0 A^{\tilde \xi} \tilde \Pi (A^1)^{-1} \tilde \Pi A^{\tilde \xi}
r_p \rangle \ne 0
\eqnlbl{intcondition}
$$
for all $\tilde \xi\in \BbbR^{d-1}\ne 0$, 
where $A^{\tilde \xi}:=\sum_{j\ne p} Df^j \xi_j$, $A^1:=Df^1$,
and $\tilde \Pi$ denotes the projection complementary to the
eigenprojection of $A^1$ onto $r_p$
(``nonresonance,'' or genuine coupling
in the sense of [M\'e.1--4,FZ,Z.5--6]).
\medskip

Here, $\epsilon>0$ as usual denotes shock strength 
$|u^\epsilon_+-u^\epsilon_-|$, 
and profiles $\bar u^\epsilon(\cdot)$ are assumed to 
converge as $\epsilon\to 0$ to $u_0$.
Since the profile problem is restricted to the normal direction,
the one-dimensional analysis of Majda\&Pego [MP]
again yields Proposition \thmref{limprofile}, with
convergence after rescaling of $\bar u^\epsilon$ to the standard
Burgers profile \eqnref{burgersprofile}.
Linearizing about $\bar u^\epsilon$,
we obtain a linearized evolution system \eqnref{ablinear}
with coefficients depending only on $x_1$;
taking the Fourier transform in directions $\tilde x:=(x_2, \dots, x_d)$,
we then obtain for each fixed $\epsilon$
a {\it family} of one-dimensional equations
$$
\hat v_t=L_{\tilde \xi} \hat v,
\eqnlbl{ablinearmult}
$$
indexed by $\tilde\xi\in \BbbR^{d-1}$.

\medskip
{\bf Definition \thmlbl{specstab}.} Following [Z.3],
we define strong spectral stability (condition (D1') of the reference)
in the multidimensional case as 
$$
\R \sigma_p (L_\xi)<0 \, \text{\rm for }\, (\tilde \xi,\lambda)\ne (0,0):
\eqnlbl{multidspecstab}
$$
equivalently, the generalized eigenvalue
equation $(L_{\tilde \xi}-\lambda)w=0$ admits no $L^2$ solutions
for $\R\lambda \ge 0$ and $(\tilde \xi,\lambda)\ne (0,0)$.
\medskip

Under hypotheses (H0)--(H4), the results of [Z.3] show that 
strong spectral spectral stability implies linearized and 
nonlinear stability with sharp rates of decay
(as discussed in [PZ], our conditions (H0)--(H4) imply conditions
(H0)--(H7) assumed in the reference).
In [PZ] we establish:

\proclaim{Theorem \thmlbl{mstab}}
Under assumptions (H0)--(H4), profiles $\bar u^\epsilon$
are strongly spectrally stable (and therefore linearly
and nonlinearly orbitally stable [ZH,Z.2])
for $\epsilon$ sufficiently small.
\endproclaim

\medskip
{\bf Remark \thmlbl{rsymm}.}
As in the inviscid small
amplitude case [M\'e.1--5,FM,Z.5--6,FZ], the condition
of simultaneous symmetrizability (H2) may be relaxed
under appropriate alternative structural conditions; however,
so far there has been identified no such conditions
of simple and general application.
Likewise, the nonresonance condition (H4)(ii) might possibly
be relaxed; however, the degenerate case that this condition
fails would appear to require a considerably more detailed analysis
than the one carried out here.
In any case, these conditions are satisfied for most systems arising
in physical applications; see [M\'e.5,Z.3,Z.5--6].
In the extreme shock case, (H4)(ii) is equivalent to a simpler
and easily verified condition identified by M\'etivier [M\'e.1--5]:
that the principal eigenvalue $a_p$ of the full symbol
$\sum_{j=1}^{d}\xi_j Df^j$ be a strictly convex, resp. concave, 
function of $\xi=(\xi_1,\dots,\xi_d)$, according as $p=1$ or $n$;
see [FZ] or (implicitly) [M\'e.1].
Note that \eqnref{intcondition} is automatically satisfied in the
extreme shock case.
The case of ``real,'' or partially parabolic, viscosity, or relaxation
may be treated similarly, following the model of the present, one-dimensional
analysis.
\bigskip

Similarly as in the one-dimensional case, Theorem \thmref{mstab}
is established by showing that stability is equivalent,
modulo pure imaginary eigenvalues as in Theorems 
\thmref{equiv} and \thmref{requiv}, to stability
of an appropriate canonical system.
However, this system is not a $2\times 2$ system as in 
the scalar multidimensional case, but
a $3\times 3$ system corresponding to the full slow
flow of an appropriate $2\times 2$ viscous conservation law,
namely, the {\it coupled Burgers-linear degenerate equations},
$$
\pmatrix u\\v\endpmatrix_t
+ \pmatrix u^2/2\\a v\endpmatrix_{x_1}
+ \pmatrix v\\u\endpmatrix_{x_2}
=
\sum B^{jk} 
\pmatrix  u \\ v 
\endpmatrix_{x_j\, x_k},
\eqnlbl{multidburgers}
$$
considered with the family of profiles
$$
\pmatrix \bar u^\epsilon\\ \bar v^\epsilon\endpmatrix (x)=
\pmatrix -\epsilon \tanh(\epsilon x/2) \\ 0 \endpmatrix,
\eqnlbl{multidprofile}
$$
where $a>0$ and $B^{jk}$ are constant, diagonal, and 
satisfy the uniform parabolicity condition
$$
\sum_{jk}\xi_j \xi_k B^{jk} \ge \theta |\xi|^2,
\quad \theta >0,
\eqnlbl{uniform}
$$
for all $\xi:=(\xi_1,\xi_2)\in \BbbR^2$,
with $B^{11}_{11}=1$.

That is, {\it in multiple dimensions, 
small amplitude behavior is not captured by any scalar model,
but generically requires a $2\times 2$ model for its expression}.
This fundamental observation was made in the inviscid setting by
M\'etivier [M\'e.1--5,FM] in his study of the weak shock limit, 
and is due to the new phenomenon of ``glancing modes'' arising in the 
multi-dimensional case;
see also [Z.5--6,FZ] for related discussions in the inviscid case.
Accordingly, the analysis of the reduced equations on the
``slow--superslow'' manifold becomes considerably
more complicated in multiple dimensions; in particular, we
can no longer rely on simple Sturm--Liouville, or maximum
principles to conclude stability of the reduced system.
The following proposition, though it does not 
(because the equations \eqnref{multidburgers} are not scale-invariant) 
directly imply stability of the reduced system, 
nonetheless indicates at an intuitive level
why the analysis can still be carried out, to yield a
full stability result and not only a reduction.
Note that the reduction to constant coefficients $a$, $A^2$,
$B^{jk}$ (achieved as in the one-dimensional
analysis by an application of Lemma \thmref{evanslimit})  
is crucial in the proof.

\proclaim{Proposition \thmlbl{modelstab}}
Shock profiles \eqnref{multidprofile} of 
any amplitude $\epsilon$ are strongly spectrally stable as solutions
of \eqnref{multidburgers}, assuming
only the partial parabolicity condition 
$$
\sum_{jk}\xi_j \xi_k B^{jk} \ge \theta |\xi|^2
\pmatrix 1 & 0\\ 0& 0\endpmatrix,
\quad \theta >0.
\eqnlbl{weakuniform}
$$
\endproclaim

\medskip
{\bf Proof.}
We use a straightforward energy estimate modifying
that used in [MN,Go.1] to treat the one-dimensional
genuinely nonlinear scalar case.  
Linearizing \eqnref{multidburgers} about a stationary profile
\eqnref{multidprofile}, we obtain the generalized eigenvalue
equation
$$
(\lambda +i\xi_2 A^2 + \xi_2^2 B^{22})w
+ (Aw)'+ i\xi_2 (B^{12}+B^{21})w'= B^{11}w'',
\eqnlbl{modelevalue}
$$
where
$$
A^1:= \pmatrix \bar u(x_1) & 0 \\
0 & a\\
\endpmatrix
\qquad
A^2:= \pmatrix 0 & 1 \\
1 & 0\\
\endpmatrix;
\eqnlbl{modelAs}
$$
here, as usual, `$'$' denotes $\partial/\partial x_1$.

Examining the symbols
$$
P_\pm(i\xi):=i\xi_1A^1_\pm+i\xi_2 A^2 +
\xi^2_1 B^{11}
+ i\xi_2 (B^{12}+B^{21})+ \xi_2^2 B^{22}
$$
of the constant coefficient limiting equations of \eqnref{modelevalue}
as $x_1\to \pm \infty$,
we find from \eqnref{weakuniform} that 
$$
\R P_\pm \ge \theta |\xi|^2
\pmatrix 1&0\\0&0 \endpmatrix,
\eqnlbl{symbol}
$$
whence, for $\xi\ne 0$,  
$P_\pm$ has no pure imaginary eigenvalues $\lambda$ unless associated
with an eigenvector of form $(0,1)^t$, in which case we find by
direct calculation that $\xi_2=0$ and $\lambda=-a\xi_1$ and 
\eqnref{modelevalue} decouples into the one-dimensional Burgers
eigenvalue equation and a trivial, constant--coefficient equation.
In either case, we find from standard considerations
(see, e.g., [He]) that $w$ and all derivatives decay exponentially 
at spatial infinity for $w\in L^2$ and $(\xi,\lambda)\ne (0,0)$.

Similarly, from
$$
\R(\lambda +i\xi_2 A^2 + \xi_2^2 B^{22})=
\R(\lambda)I + \xi_2^2 \R(B^{22}) \ge \theta_2 \xi_2^2
\pmatrix 1&0\\0&0\endpmatrix
$$
for $\theta_2>0$, by \eqnref{uniform}, we easily find that
$(\lambda +i\xi_2 A^2 + \xi_2^2 B^{22})$ is invertible
for all real $\xi_2$ and $\R \lambda \ge 0$ such that
$(\xi_2, \lambda)\ne (0,0)$, whence 
$\int_{-\infty}^{+\infty}w(z)dz=0$; for, this implies that the
kernel, if it exists, is spanned by $(0,1)^t)$, 
yielding $\xi_2=\lambda=0$ by direct calculation.
Combining these observations, and integrating \eqnref{modelevalue} from
$-\infty$ to $+\infty$, 
we obtain 
$$
(\lambda +i\xi_2 A^2 + \xi_2^2 B^{22})
\int_{-\infty}^{+\infty}w(z)dz=0
\eqnlbl{zeromass}
$$
and therefore $ \int_{-\infty}^{+\infty}w(z)dz=0$.

Thus, we may introduce the integrated variable 
$W(x_1):=\int_{-\infty}^{x_1} w(z)dz $, with $W\in L^2$
satisfying the integrated eigenvalue equation
$$
(\lambda +i\xi_2 A^2 + \xi_2^2 B^{22})W
+ AW'+ i\xi_2 (B^{12}+B^{21})W'= B^{11}W'',
\eqnlbl{modelintevalue}
$$
with $W$ and derivatives again decaying exponentially at spatial infinity.
Taking the real part of the complex inner product of $W$ against
\eqnref{modelintevalue}, integrating by parts, and using \eqnref{uniform},
we obtain the estimate
$$
\aligned
\R \lambda |W|^2
-\int \bar u' |W1|^2 dx&= 
-\langle W',B^{11} W'\rangle
+\langle i\xi_2 W, B^{12} W'\rangle\\
&\quad+\langle W', B^{12} i\xi_2 W\rangle
+\langle i\xi_2 W, B^{22} i\xi_2 W \rangle
\\
&\le -\theta(|W_1'|^2+ \xi_2^2|W_1|^2),
\endaligned
$$
$\theta>0$, a contradiction of the fact that $\bar u'<0$,
unless $|W_1'|\equiv \xi_2 |W_1|\equiv 0$, hence $|W_1|\equiv 0$, 
and $\R \lambda = 0$.
But, this case may be disposed of by observing that \eqnref{modelevalue}
then reduces to a single equation in $W_2$,
which is constant-coefficient and therefore
supports no $L^2$ solutions.
\qed

\bigskip
\Refs
\medskip\noindent
[AGJ] J. Alexander, R. Gardner and C.K.R.T. Jones,
{\it A topological invariant arising in the analysis of
traveling waves}, J. Reine Angew. Math. 410 (1990) 167--212.
\medskip\noindent
[Ce] C. Cercignani,
{\it The Boltzmann equation and its applications,}
Applied Mathematical Sciences, 67. Springer-Verlag, New York (1988)
xii+455 pp. ISBN: 0-387-96637-4.
\medskip\noindent
[CLL] G.-Q. Chen, D.C. Levermore, and T.-P. Liu, 
{\it Hyperbolic conservation laws with 
stiff relaxation terms and entropy,} Comm. Pure
Appl. Math. 47 (1994), no. 6, 787--830. 
\medskip\noindent
[CL] E.A. Coddington and N. Levinson, {\it Theory of ordinary
differential equations,} McGraw-Hill Book Company, Inc., 
New York-Toronto-London (1955) xii+429 pp.
\medskip\noindent
[Co] W. A. Coppel,
{\it Stability and asymptotic behavior of differential equations},
D.C. Heath and Co., Boston, MA (1965) viii+166 pp. 
\medskip\noindent
[Dre] W. Dreyer, 
{\it Maximisation of the entropy in nonequilibrium,}
J. Phys. A 20 (1987), no. 18, 6505--6517. 
\medskip\noindent
[E.1] J.W. Evans,
{\it Nerve axon equations: I. Linear approximations,}
Ind. Univ. Math. J. 21 (1972) 877--885.
\medskip\noindent
[E.2] J.W. Evans,
{\it Nerve axon equations: II. Stability at rest,}
Ind. Univ. Math. J. 22 (1972) 75--90.
\medskip\noindent
[E.3] J.W. Evans,
{\it Nerve axon equations: III. Stability of the nerve impulse,}
Ind. Univ. Math. J. 22 (1972) 577--593.
\medskip\noindent
[E.4] J.W. Evans,
{\it Nerve axon equations: IV. The stable and the unstable impulse,}
Ind. Univ. Math. J. 24 (1975) 1169--1190.
\medskip\noindent
[FM] J. Francheteau and G. M\'etivier,
{\it Existence de chocs faibles pour des syst\'emes quasi-lin\'eaires
hyperboliques multidimensionnels,}
C.R.AC.Sc> Paris, 327 S\'erie I (1998) 725--728.
\medskip\noindent
[FreL] H. Freist\"uhler and T.-P. Liu, 
{\it Nonlinear stability of overcompressive shock waves 
in a rotationally invariant system of viscous conservation laws,}
Comm. Math. Phys. 153 (1993) 147--158.
\medskip\noindent
[FreS] H. Freist\"uhler and P. Szmolyan, 
{\it Spectral stability of small shock waves, I,}
preprint (received March 30, 2002).
\medskip\noindent
[FZe] H. Freist\"uhler and Y. Zeng,
{\it Shock profiles for systems of balance
laws with relaxation,} preprint (1998).
\medskip\noindent
[FZ] H. Freist\"uhler and K. Zumbrun,
{\it Long-wave stability of small-amplitude 
overcompressive shock profiles, with applications to MHD,}
in preparation.
\medskip\noindent
[GJ.1] R. Gardner and C.K.R.T. Jones,
{\it A stability index for steady state solutions of
boundary value problems for parabolic systems},
J. Diff. Eqs. 91 (1991), no. 2, 181--203. 
\medskip\noindent
[GJ.2] R. Gardner and C.K.R.T. Jones,
{\it Traveling waves of a perturbed diffusion equation
arising in a phase field model}, 
Ind. Univ. Math. J. 38 (1989), no. 4, 1197--1222.
\medskip\noindent
[GJ.3] R. Gardner and C.K.R.T. Jones,
{\it Stability of one-dimensional waves in weak and singular limits,}
Viscous profiles and numerical methods for shock waves 
(Raleigh, NC, 1990), 32--48, SIAM, Philadelphia, PA, 1991. 
\medskip\noindent
[GZ] R. Gardner and K. Zumbrun,
{\it The Gap Lemma and geometric criteria for instability
of viscous shock profiles}, 
Comm. Pure Appl.  Math. 51 (1998), no. 7, 797--855. 
\medskip\noindent
[G] P. Godillon,
{\it Linear stability of shock profiles for systems of conservation laws 
with semi-linear relaxation,} (English. English summary) 
Phys. D 148 (2001), no. 3-4, 289--316. 
\medskip\noindent
[Go.1] J. Goodman, {\it Nonlinear asymptotic stability of viscous
shock profiles for conservation laws,}
Arch. Rational Mech. Anal. 95 (1986), no. 4, 325--344.
\medskip\noindent
[Go.2] J. Goodman, 
{\it Remarks on the stability of viscous shock waves}, in:
Viscous profiles and numerical methods for shock waves 
(Raleigh, NC, 1990), 66--72, SIAM, Philadelphia, PA, (1991).
\medskip\noindent
[He] D. Henry,
{\it Geometric theory of semilinear parabolic equations},
Lecture Notes in Mathematics, Springer--Verlag, Berlin (1981),
iv + 348 pp.
\medskip\noindent
[HoZ.1] D. Hoff and K. Zumbrun, 
{\it Stability and asymptotic behavior of multi-dimensional 
scalar viscous shock fronts,}  
Indiana Univ. Math. J. 49 (2000) 427--474. 
\medskip\noindent
[HoZ.2] D. Hoff and K. Zumbrun, 
{\it Pointwise Green's function bounds for multi-dimensional 
scalar viscous shock fronts,}  to appear, J. Diff. Eq. (2002).
\medskip\noindent
[H.1] P. Howard,
{\it Pointwise estimates on the Green's function 
for a scalar linear convection-diffusion equation,}
J. Differential Equations 155 (1999), no. 2, 327--367.
\medskip\noindent
[H.2] P. Howard,
{\it Pointwise methods for stability of a scalar conservation law,}
Doctoral thesis (1998).
\medskip\noindent
[H.3] P. Howard,
{\it Pointwise Green's function approach to 
stability for scalar conservation laws,}
Comm. Pure Appl. Math. 52 (1999), no. 10, 1295--1313.
\medskip\noindent
[HZ] P. Howard and K. Zumbrun,
{\it Pointwise estimates for dispersive-diffusive shock waves,}
to appear, Arch. Rational Mech. Anal. 
\medskip\noindent
[Hu] J. Humpherys, 
{\it Stability of Jin--Xin relaxation shocks,} preprint (2000).
\medskip\noindent
[HuZ], J. Humpherys and K. Zumbrun, 
{\it Spectral stability of small-amplitude shock profiles for
dissipative hyperbolic-parabolic systems of conservation laws,} 
Z. angew. Math. Phys. 53 (2002) 20--34.
\medskip\noindent
[JX] S. Jin and Z. Xin, 
{\it The relaxation schemes for systems of conservation laws 
in arbitrary space dimensions,}
Comm. Pure Appl.  Math. 48 (1995), no. 3, 235--276.
\medskip\noindent
[J] C.K.R.T. Jones,
{\it Stability of the travelling wave solution of the FitzHugh--Nagumo system},
 Trans. Amer. Math. Soc.  286 (1984), no. 2, 431--469.
\medskip\noindent
[JGK] C. K. R. T. Jones, R. A. Gardner, and T. Kapitula,
{\it Stability of travelling waves
for non-convex scalar viscous conservation laws},
Comm. Pure Appl. Math. 46 (1993) 505--526.
\medskip\noindent
[K.1] T. Kapitula,
{\it Stability of weak shocks in $\lambda$--$\omega$ systems},
Indiana Univ. Math. J. 40 (1991), no. 4, 1193--12.
\medskip\noindent
[K.2] T. Kapitula,
{\it On the stability of travelling waves in weighted $L^\infty$ spaces},
J. Diff. Eqs. 112 (1994), no. 1, 179--215.
\medskip\noindent
[KS] T. Kapitula and B. Sandstede,
{\it Stability of bright solitary-wave solutions 
to perturbed nonlinear Schrödinger equations,} Phys. D 124
(1998), no. 1-3, 58--103.
\medskip\noindent
[Kaw] S. Kawashima,
{\it Systems of a hyperbolic--parabolic composite type,
with applications to the equations of magnetohydrodynamics},
thesis, Kyoto University (1983).
\medskip\noindent
[KM] S. Kawashima and A. Matsumura,
{\it Asymptotic stability of traveling wave solutions of systems 
for one-dimensional gas motion,} Comm.
Math. Phys. 101 (1985), no. 1, 97--127.
\medskip\noindent
[KMN] S. Kawashima, A. Matsumura, and K. Nishihara,
{\it Asymptotic behavior of solutions for the equations of a 
viscous heat-conductive gas,}
Proc.  Japan Acad. Ser. A Math. Sci. 62 (1986), no. 7, 249--252. 
\medskip\noindent
[Kat] T. Kato,
{\it Perturbation theory for linear operators},
Springer--Verlag, Berlin Heidelberg (1985).
\medskip\noindent
[Lev]
D.C. Levermore, 
{\it Moment closure hierarchies for kinetic theories,}
(English. English summary) 
J. Statist. Phys. 83 (1996), no. 5-6, 1021--1065. 
\medskip\noindent
[Liu] H. Liu,
{\it Asymptotic stability of relaxation shock
profiles for hyperbolic conservation laws,}
preprint (2000).  
\medskip\noindent
[L.1] T.-P. Liu,
{\it Nonlinear stability of shock waves for viscous conservation laws},
Mem. Amer. Math. Soc. 56 (1985), no. 328, v+108 pp.
\medskip\noindent
[L.2] T.-P. Liu,
{\it Pointwise convergence to shock waves for viscous conservation laws},
Comm. Pure Appl. Math. 50 (1997), no. 11, 1113--1182.
\medskip\noindent
[L.3] T.-P. Liu, 
{\it Hyperbolic conservation laws with relaxation,}
Comm. Math. Phys. 108 (1987), no. 1, 153--175. 
[M.1], A. Majda,
{\it The stability of multi-dimensional shock fronts -- a
new problem for linear hyperbolic equations,}
Mem. Amer. Math. Soc. 275 (1983).
\medskip\noindent
[M.2], A. Majda,
{\it The existence of multi-dimensional shock fronts,}
Mem. Amer. Math. Soc. 281 (1983).
\medskip\noindent
[M.3] A. Majda,
{\it Compressible fluid flow and systems of conservation laws in several
space variables,} Springer-Verlag, New York (1984), viii+ 159 pp.
\medskip\noindent
[MP] A. Majda and R. Pego,
{\it Stable viscosity matrices for systems of conservation laws},
J. Diff. Eqs. 56 (1985) 229--262.
\medskip\noindent
[MaN] C. Mascia and R. Natalini, 
{\it $L\sp 1$ nonlinear stability of traveling waves for a 
hyperbolic system with relaxation,} J. Differential
Equations 132 (1996), no. 2, 275--292.
\medskip\noindent
[MZ.1] C. Mascia and K. Zumbrun,  Relaxation shocks, 
{\it Pointwise Green's function bounds and stability of relaxation shocks,}
preprint (2001).
\medskip\noindent
[MZ.2] C. Mascia and K. Zumbrun,  
{\it Stability of viscous shock profiles for dissipative 
hyperbolic-parabolic systems of conservation laws,}
preprint (2001).
\medskip\noindent
[MN] A. Matsumura and K. Nishihara, 
{\it On the stability of travelling wave solutions of a one-dimensional 
model system for compressible viscous gas,}
Japan J. Appl. Math. 2 (1985), no. 1, 17--25.
\medskip\noindent
[M\'e.1] G. M\'etivier, {\it Stability of multidimensional weak shocks,}
Comm. Partial Diff. Equ. 15 (1990) 983--1028.
\medskip\noindent
[M\'e.2] G. M\'etivier, {\it Interaction de deux choc pour un
syst\'eme de deux lois de conservation, en dimension deux d'espace,}
Trans. Amer. Math. Soc. 296 (1986) 431--479.
\medskip\noindent
[M\'e.3] G. M\'etivier, {\it Ondes soniques,}
J. Math. pures et appl. 70 (1991) 197--268.
\medskip\noindent
[M\'e.4] G. M\'etivier,
{\it The block structure condition for symmetric hyperbolic systems,}
preprint (1999).
\medskip\noindent
[M\'e.5] G. M\'etivier,
{\it Stability of multidimensional shocks,}
Advances in the theory of shock waves, 25--103, 
Progr. Nonlinear Differential Equations Appl., 47, 
Birkhäuser Boston, Boston, MA, 2001.
\medskip\noindent
[M\'eZ], G. M\'etivier and K Zumbrun,
{\it Large Viscous Boundary Layers for Noncharacteristic
Nonlinear Hyperbolic Problems,}
preprint (2002).
\medskip\noindent
[MR] I. M\"uller and T. Ruggeri,
{\it Rational extended thermodynamics,}
Second edition, with supplementary chapters by H. Struchtrup 
and Wolf Weiss, Springer Tracts in Natural Philosophy, 37, 
Springer-Verlag, New York (1998) xvi+396 pp. ISBN: 0-387-98373-2.
\medskip\noindent
[N] R. Natalini, 
{\it Recent mathematical results on hyperbolic relaxation problems,}
TMR Lecture Notes (1998).
\medskip\noindent
[OZ.1] M. Oh and K. Zumbrun,
{\it Stability of periodic solutions of viscous conservation laws:
Analysis of the Evans function,} to appear, Archive for Rat. Mech.
\medskip\noindent
[OZ.2] M. Oh and K. Zumbrun,
{\it Stability of periodic solutions of viscous conservation laws:
Pointwise bounds on the Green function,} 
to appear, Archive for Rat. Mech.
\medskip\noindent
[P] R. Pego, 
{\it Stable viscosities and shock profiles for 
systems of conservation laws,}
Trans. Amer. Math. Soc. 282 (1984) 749--763.
\medskip\noindent
[PI] T.  Platkowski and R. Illner, 
{\it Discrete velocity models of the Boltzmann equation: 
a survey on the mathematical aspects of the
theory,} SIAM Rev. 30 (1988), no. 2, 213--255. 
\medskip\noindent
[Pl] R. Plaza,
{\it Center manifold and spectral stability of weak viscous
shock waves,} in preparation (2002).
\medskip\noindent
[PZ] R. Plaza and K. Zumbrun,
{\it Multidimensional stability of small amplitude viscous
shock fronts,} in preparation. 
\medskip\noindent
[Sat] D. Sattinger,
{\it On the stability of waves of nonlinear parabolic systems},
Adv. Math. 22 (1976) 312--355.
\medskip\noindent
[Sm] J. Smoller, 
{\it Shock waves and reaction--diffusion equations,} 
Second edition, Grundlehren der Mathematischen Wissenschaften
[Fundamental Principles of Mathematical Sciences], 258. Springer-Verlag, 
New York, 1994. xxiv+632 pp. ISBN: 0-387-94259-9. 
\medskip\noindent
[S] A.  Szepessy, 
{\it On the stability of Broadwell shocks,}
Nonlinear evolutionary partial differential equations 
(Beijing, 1993), 403--412, 
AMS/IP Stud. Adv. Math., 3, Amer. Math. Soc., Providence, RI, 1997. 
\medskip\noindent
[SX.1] A. Szepessy and Z. Xin,
{\it Nonlinear stability of viscous shock waves,}
Arch. Rat. Mech. Anal. 122 (1993) 53--103.
\medskip\noindent
[SX.2] A. Szepessy and Z. Xin,
{\it Stability of Broadwell shocks,}
unpublished manuscript.
\medskip\noindent
[Wh] G.B. Whitham, 
{\it Linear and nonlinear waves,}
Reprint of the 1974 original, Pure and Applied Mathematics,
A Wiley-Interscience Publication, John Wiley \& Sons, Inc., 
New York, 1999. xviii+636 pp. ISBN: 0-471-35942-4. 
\medskip\noindent
[Yo.4] W.-A. Yong, 
{\it Basic aspects of hyperbolic relaxation systems,}
Advances in the theory of shock waves, 259--305,
Progr. Nonlinear Differential Equations Appl., 47, 
Birkhäuser Boston, Boston, MA, 2001.
\medskip\noindent
[YoZ] W.-A. Yong and K. Zumbrun,
{\it Existence of relaxation shock profiles 
for hyperbolic conservation laws,}
SIAM J. Appl. Math. 60 (2000), no. 5, 1565--1575 (electronic).
\medskip\noindent
[Y] K. Yosida, {\it Functional analysis,} 
Reprint of the sixth (1980) edition, Classics in Mathematics, 
Springer-Verlag, Berlin, 1995, xii+501 pp. ISBN: 3-540-58654-7.
\medskip\noindent
[Ze] Y. Zeng,
{\it Gas dynamics in thermal nonequilibrium 
and general hyperbolic systems with relaxation,}
Arch. Ration. Mech. Anal.  150 (1999), no. 3, 225--279.
\medskip\noindent
[ZH] K. Zumbrun and P. Howard,
{\it Pointwise semigroup methods and stability of viscous shock waves,}
Indiana Mathematics Journal V47 (1998) no. 4, 741--871.
\medskip\noindent
[Z.1] K. Zumbrun,  {\it Stability of viscous shock waves},
Lecture Notes, Indiana University (1998).
\medskip\noindent
[Z.2] K. Zumbrun, {\it Refined Wave--tracking and Nonlinear 
Stability of Viscous Lax Shocks}, Methods Appl. Anal. (2000)
747--768.
\medskip\noindent
[Z.3] K. Zumbrun, {\it Multidimensional stability of
planar viscous shock waves}, 
Birkhauser's Series: Progress in Nonlinear Differential
Equations and their Applications (2001), 207 pp.
\medskip\noindent
[Z.4] K. Zumbrun, 
{\it Dynamical stability of phase transitions in the $p$-system 
with viscosity-capillarity,}
SIAM J. Appl. Math. 60 (2000), 1913--1924.
\medskip\noindent
[Z.5] K. Zumbrun, {\it Multidimensional stability of
planar viscous shock fronts,} 
Lecture Notes, Indiana University (2000).
\medskip\noindent
[Z.6] K. Zumbrun, {\it Multidimensional stability of 
planar viscous shock fronts of compressible Navier--Stokes equations
and related equations of compressible flow,} 
in {\it Handbook of fluid Dynamics}, in preparation.
\medskip\noindent
[ZS] K. Zumbrun and D. Serre,
{\it Viscous and inviscid stability of multidimensional 
planar shock fronts,} Indiana Univ. Math. J. 48 (1999), no. 3,
937--992.
\endRefs
\enddocument

%% file: eqnumams.tex
\newif\ifproofmode                      
\proofmodefalse				

\newif\ifforwardreference		
\forwardreferencefalse			

\newif\ifchapternumbers			
\chapternumbersfalse			

\newif\ifcontinuousnumbering		
\continuousnumberingfalse		

\newif\iffigurechapternumbers		
\figurechapternumbersfalse		

\newif\ifcontinuousfigurenumbering	
\continuousfigurenumberingfalse		

\font\eqsixrm=cmr6			
\def\marginstyle{\eqsixrm}		

\newtoks\chapletter			
\newcount\chapno			
\newcount\eqlabelno			
\newcount\figureno			

\chapno=0
\eqlabelno=0
\figureno=0


\def\chapfolio{\ifnum \chapno>0 \the\chapno \else \the\chapletter \fi}


\def\bumpchapno{\ifnum \chapno>-1 \global \advance \chapno by 1
	\else \global \advance \chapno by -1 \setletter\chapno \fi
	\ifcontinuousnumbering \else \global\eqlabelno=0 \fi
	\ifcontinuousfigurenumbering \else \global\figureno=0 \fi}

%


%

\def\tempsetletter#1{\ifcase-#1 {}\or{} \or\chapletter={A}\or\chapletter={B}
  \or\chapletter={C} \or\chapletter={D} \or\chapletter={E}
  \or\chapletter={F} \or\chapletter={G} \or\chapletter={H}
  \or\chapletter={I} \or\chapletter={J} \or\chapletter={K}
  \or\chapletter={L} \or\chapletter={M} \or\chapletter={N}
  \or\chapletter={O} \or\chapletter={P} \or\chapletter={Q}
  \or\chapletter={R} \or\chapletter={S} \or\chapletter={T}
  \or\chapletter={U} \or\chapletter={V} \or\chapletter={W}
  \or\chapletter={X} \or\chapletter={Y} \or\chapletter={Z}\fi}

%

\def\chapshow#1{\ifnum #1>0 \relax #1%
   \else {\tempsetletter{\number#1}\chapno=#1 \chapfolio} \fi}

%
\def\today{\ifcase\month\or
January\or February\or March\or April\or May\or June\or
July\or August\or September\or October\or November\or December\fi
\space\number\day, \number\year}

%

\def\initialeqmacro{\ifproofmode
 \headline{\tenrm \today\hfill \jobname\ --- draft\hfill\folio}
     \hoffset=-1cm \immediate\openout2=allcrossreferfile \fi
 \ifforwardreference \input labelfile
     \ifproofmode \immediate\openout1=labelfile \fi \fi}


%

\def\chaplabel#1{\bumpchapno \ifproofmode \ifforwardreference
   \immediate\write1{\noexpand\expandafter\noexpand\def
   \noexpand\csname CHAPLABEL#1\endcsname{\the\chapno}}\fi\fi
   \global\expandafter\edef\csname CHAPLABEL#1\endcsname
   {\the\chapno}\ifproofmode\llap{\hbox{\marginstyle #1\ }}\fi\chapfolio}

%
\def\eqnum{\global\advance\eqlabelno by 1
   \eqno(\ifchapternumbers\chapfolio.\fi\the\eqlabelno)}

\def\eqlabel#1{\global\advance\eqlabelno by 1 \ifproofmode\ifforwardreference
 \immediate\write1{\noexpand\expandafter\noexpand\def
 \noexpand\csname EQLABEL#1\endcsname{\the\chapno.\the\eqlabelno?}}\fi\fi
 \global\expandafter\edef\csname EQLABEL#1\endcsname
 {\the\chapno.\the\eqlabelno?} \eqno(\ifchapternumbers\chapfolio.\fi
 \the\eqlabelno)\ifproofmode\rlap{\hbox{\marginstyle #1}}\fi}

\def\leqlabel#1{\global\advance\eqlabelno by 1 \ifproofmode\ifforwardreference
 \immediate\write1{\noexpand\expandafter\noexpand\def
 \noexpand\csname EQLABEL#1\endcsname{\the\chapno.\the\eqlabelno?}}\fi\fi
 \global\expandafter\edef\csname EQLABEL#1\endcsname
 {\the\chapno.\the\eqlabelno?} \leqno(\ifchapternumbers\chapfolio.\fi
 \the\eqlabelno)\ifproofmode\rlap{\hbox{\marginstyle #1}}\fi}

\def\eqalignnum{\global\advance\eqlabelno by 1
   &(\ifchapternumbers\chapfolio.\fi\the\eqlabelno)}

\def\eqalignlabel#1{\global\advance\eqlabelno by1 \ifproofmode
 \ifforwardreference\immediate\write1{\noexpand\expandafter\noexpand\def
 \noexpand\csname EQLABEL#1\endcsname
     {\the\chapno.\the\eqlabelno?}}\fi\fi
 \global\expandafter\edef\csname EQLABEL#1\endcsname
 {\the\chapno.\the\eqlabelno?}\ifchapternumbers\chapfolio.\fi
 \the\eqlabelno\ifproofmode\rlap{\hbox{\marginstyle
 #1}}\fi}

\def\eqref#1{(\ifundefined{EQLABEL#1}***\ifproofmode\ifforwardreference)%
   \else \write16{ ***Undefined Equation Reference #1*** }\fi
   \else \write16{ ***Undefined Equation Reference #1*** }\fi
   \else \edef\LABxx{\getlabel{EQLABEL#1}}%
   \def\LAByy{\expandafter\stripchap\LABxx}\ifchapternumbers
   \chapshow{\LAByy}.\expandafter\stripeq\LABxx
   \else\ifnum \number\LAByy=\chapno \relax\expandafter\stripeq\LABxx
   \else\chapshow{\LAByy}.\expandafter\stripeq\LABxx\fi\fi)\fi
   \ifproofmode\write2{Equation #1}\fi}

%

\def\fignum{\global\advance\figureno by 1 \relax
   \iffigurechapternumbers\chapfolio.\fi\the\figureno}\

\def\figlabel#1{\global\advance\figureno by 1\relax
 \ifproofmode\ifforwardreference
 \immediate\write1{\noexpand\expandafter\noexpand\def
 \noexpand\csname FIGLABEL#1\endcsname{\the\chapno.\the\figureno?}}\fi\fi
 \global\expandafter\edef\csname FIGLABEL#1\endcsname
 {\the\chapno.\the\figureno?}\iffigurechapternumbers\chapfolio.\fi
 \ifproofmode$^{\hbox{\marginstyle #1}}$\relax\fi\the\figureno}

\def\figref#1{\ifundefined{FIGLABEL#1}!!!!\ifproofmode\ifforwardreference)%
   \else \write16{ ***Undefined Equation Reference #1*** }\fi
   \else \write16{ ***Undefined Equation Reference #1*** }\fi
   \else \edef\LABxx{\getlabel{FIGLABEL#1}}%
   \def\LAByy{\expandafter\stripchap\LABxx}%
   \iffigurechapternumbers\chapshow{\LAByy}.\expandafter\stripeq\LABxx
   \else\ifnum\number\LAByy=\chapno \relax\expandafter\stripeq\LABxx
   \else\chapshow{\LAByy}.\expandafter\stripeq\LABxx\fi\fi
   \ifproofmode\write2{Figure #1}\fi\fi}

%

%

\def\getlabel#1{\csname#1\endcsname}
\def\ifundefined#1{\expandafter\ifx\csname#1\endcsname\relax}
\def\stripchap#1.#2?{#1}
\def\stripeq#1.#2?{#2}

\figurechapternumberstrue  

\chapternumberstrue        

\def\thmlbl#1{\figlabel{#1}}
\def\thmref{\figref}
\def\eqnlbl#1{\leqlabel{#1}}

\def\eqnref#1{\eqref{#1}}
\def\sectionnumber{\chapno}
\def\theoremnumber{\figureno}
\def\equationnumber{\eqlabelno}